%%%%%%%%%%%%%%%%%%%%%%%%%%%%%%%%%%%%%%%%%%%%%%%%%%%%%%%%%%%%%%%%%%%%%%%%%%%%%
%%%%%%%%%%%%%%%%%%%%%%%%%%%%%%%%%%%%%%%%%%%%%%%%%%%%%%%%%%%%%%%%%%%%%%%%%%%%%
%%%%%%%%%%%%%%%%%%%%%%%%%%%%%%%%%%%%%%%%%%%%%%%%%%%%%%%%%%%%%%%%%%%%%%%%%%%%%
%%%%%%%%%%%%  UNITARY SPHERICAL HIGHEST WEIGHT REPRESENTATIONS   %%%%%%%%%%%%
%%%%%%%%%%%%                       BY                            %%%%%%%%%%%%
%%%%%%%%%%%%       BERNHARD KROETZ AND KARL-HERMANN NEEB         %%%%%%%%%%%%
%%%%%%%%%%%%%%%%%%%%%%%%%%%%%%%%%%%%%%%%%%%%%%%%%%%%%%%%%%%%%%%%%%%%%%%%%%%%%
%%%%%%%%%%%%%%%%%%%%%%%%%%%%%%%%%%%%%%%%%%%%%%%%%%%%%%%%%%%%%%%%%%%%%%%%%%%%%
%%%%%%%%%%%%%%%%%%%%%%%%%%%%%%%%%%%%%%%%%%%%%%%%%%%%%%%%%%%%%%%%%%%%%%%%%%%%%

%VERSION OF august 24, 1995
\input amssym.def
\input amssym.tex
% Achtung: Bei Verwendung von jltmac werden die Fontfamilien masm und
% msbm mehrfach aufgerufen. Das kann zu Kapazitaetsschwierigkeiten
% fuehren. KH 27.2.95

%\def\lit{\sectionheadline{\bf References}
%\frenchspacing
%\entries\par}

\def\item#1{\vskip1.3pt\hang\textindent {\rm #1}}% THIS REPLACES KNUTH'S DEF'N

\def\itemitem#1{\vskip1.3pt\indent\hangindent2\parindent\textindent {\rm #1}}
                                  % THIS REPLACES KNUTH'S DEF'N

\tolerance=300
\pretolerance=200
\hfuzz=1pt
\vfuzz=1pt

% Print out with -x300 -y50

% Offsetwerte fuer Ausdrucke in Erlangen, hoffset um 0.6 in groesser machen
\hoffset=0.6in
\voffset=0.8in

%\baselineskip                 
\hsize=5.8 true in 

%%In Deutschland \vsize=9.2 %%

\vsize=8.5 true in
%\baselineskip
\parindent=25pt
\mathsurround=1pt
\parskip=1pt plus .25pt minus .25pt
\normallineskiplimit=.99pt

\countdef\revised=100
\mathchardef\emptyset="001F % THIS REPLACES KNUTH'S DEFINITION
\chardef\ss="19
\def\3{\ss}
\def\anf{$\lower1.2ex\hbox{"}$}
\def\frac#1#2{{#1 \over #2}}
\def\>{>\!\!>}
\def\<{<\!\!<}

\def\ssarr{\hbox to 30pt{\rightarrowfill}}
\def\sarr{\hbox to 40pt{\rightarrowfill}}
\def\arr{\hbox to 60pt{\rightarrowfill}}
\def\larr{\hbox to 60pt{\leftarrowfill}}
\def\Arr{\hbox to 80pt{\rightarrowfill}}

{}

\def\ad{\mathop{\rm ad}\nolimits}

\def\Aut{\mathop{\rm Aut}\nolimits}

\def\der{\mathop{\rm der}\nolimits}
\def\det{\mathop{\rm det}\nolimits}

\def\End{\mathop{\rm End}\nolimits}

\def\GL{\mathop{\rm GL}\nolimits}

\def\Herm{\mathop{\rm Herm}\nolimits}

\def\Hol{\mathop{\rm Hol}\nolimits}%
\def\Hom{\mathop{\rm Hom}\nolimits}%
 % USED FOR IDENTITY FUNCTION

% USED FOR IMAGINARY PART OF COMPLEX NUMBERS

\def\Pol{\mathop{\rm Pol}\nolimits}

% USED FOR REAL PART OF COMPLEX NUMBERS

%\def\reach{\mathop{\rm reach}\nolimits}
%\def\SAut{\mathop{\rm SAut}\nolimits}
%\def\SEnd{\mathop{\rm SEnd}\nolimits}

\def\SO{\mathop{\rm SO}\nolimits}
\def\span{\mathop{\rm span}\nolimits}

\def\Sp{\mathop{\rm Sp}\nolimits}
\def\Spec{\mathop{\rm Spec}\nolimits}

\def\supp{\mathop{\rm supp}\nolimits}
\def\tr{\mathop{\rm tr}\nolimits}% USED FOR TRACE OF MATRIX

\def\0{{\bf 0}}
\def\1{{\bf 1}}

\def\a{{\frak a}}

\def\b{{\frak b}}

\def\e{{\frak e}}

\def\g{{\frak g}}
\def\gl{{\frak {gl}}}
\def\h{{\frak h}}

\def\k{{\frak k}}

\def\n{{\frak n}}

\def\p{{\frak p}}
\def\q{{\frak q}}

\def\sp{{\frak {sp}}}

\def\su{{\frak {su}}}
\def\so{{\frak {so}}}
\def\sL{{\frak {sl}}}
\def\t{{\frak t}}

\def\z{{\frak z}}

\def\C{{\Bbb C}}

\def\N{{\Bbb N}}

\def\R{{\Bbb R}} 
 
\def\Z{{\Bbb Z}} 

\def\:{\colon}  %8.5.92
\def\.{{\cdot}}
\def\|{\Vert}
\def\bsk{\bigskip}

\def\giantskip{\vskip2\bigskipamount}
\def\gsk{\giantskip}
\def \la {\langle}
\def\msk{\medskip}
\def \ra {\rangle}
\def \res {\!\mid\!\!}

\def\bbr{\bigbreak}
\def\giantbreak{\par \ifdim\lastskip<2\bigskipamount \removelastskip
         \penalty-400 \giantskip\fi}

\def\nin{\noindent}
\def\cen{\centerline}
\def\pagebreak{\vskip 0pt plus 0.0001fil\break}
\def\linebreak{\break}

\def\hat{\widehat}

\def\derat#1{{d \over dt} \hbox{\vrule width0.5pt 
                height 5mm depth 3mm${{}\atop{{}\atop{\scriptstyle t=#1}}}$}}

\def\eps{\varepsilon}
\def\epsilon{\varepsilon}

\def\nin{\noindent}
\def\oline{\overline}

\def\pder#1,#2,#3 { {\partial #1 \over \partial #2}(#3)}
\def\pde#1,#2 { {\partial #1 \over \partial #2}}
\def\phi{\varphi}

% Besser \ltimes und \rtimes aus dem AMS-Symbols verwenden. 
%\def\sdir#1{\hbox{$\mathrel\times{\hskip -4.6pt 
%            {\vrule height 4.7 pt depth .5 pt}}\hskip 2pt_{#1}$}}

\def\subeq{\subseteq}
\def\supeq{\supseteq}

\def\tilde{\widetilde}

\font\eightrm=cmr8

% SANS SERIF 10 POINT
 %SANS SERIF 10 POINT ITALIC

\font\smc=cmcsc10
%\font\smc8=cmcsc8 
 %SLANTED TYPEWRITER 10 POINT
 %BOLD FACE MATH SYMBOLS 10 POINT
 %DUNHILL STYLE 10 POINT
 %SAN SERIF BOLD EXTENDED 10 POINT
 %USED FOR TITLES
 %USED FOR TITLES
\font\bfone=cmbx10 scaled\magstep1 %BOLDFACE AT MAGSTEP 1
\font\bftwo=cmbx10 scaled\magstep2 %BOLDFACE AT MAGSTEP 2
 %BOLDFACE AT MAGSTEP 3

\def\qed{{\unskip\nobreak\hfil\penalty50\hskip .001pt \hbox{}\nobreak\hfil
          \vrule height 1.2ex width 1.1ex depth -.1ex
           \parfillskip=0pt\finalhyphendemerits=0\medbreak}\rm}
%This is the end-of-proof sign. 
%Not to be used in display mode. 
%If you want to conclude a proof 
%at the end of a line is display mode use 

\def\qeddis{\eqno{\vrule height 1.2ex width 1.1ex depth -.1ex} $$
                   \medbreak\rm}
%BUT OMIT $$---the macro will write that

\def\Lemma #1. {\bigbreak\vskip-\parskip\noindent{\bf Lemma #1.}\quad\it}

\def\Sublemma #1. {\bigbreak\vskip-\parskip\noindent{\bf Sublemma #1.}\quad\it}

\def\Proposition #1. {\bigbreak\vskip-\parskip\noindent{\bf Proposition #1.}
\quad\it}

\def\Corollary #1. {\bigbreak\vskip-\parskip\nin{\bf Corollary #1.}
\quad\it}

\def\Theorem #1. {\bigbreak\vskip-\parskip\noindent{\bf Theorem #1.}
\quad\it}

\def\Definition #1. {\rm\bigbreak\vskip-\parskip\noindent{\bf Definition #1.}
\quad}

\def\Remark #1. {\rm\bigbreak\vskip-\parskip\noindent{\bf Remark #1.}\quad}

\def\Example #1. {\rm\bigbreak\vskip-\parskip\noindent{\bf Example #1.}\quad}

\def\Problems #1. {\bigbreak\vskip-\parskip\noindent{\bf Problems #1.}\quad}
\def\Problem #1. {\bigbreak\vskip-\parskip\noindent{\bf Problems #1.}\quad}

\def\Conjecture #1. {\bigbreak\vskip-\parskip\noindent{\bf Conjecture #1.}\quad}

\def\Proof#1.{\rm\par\ifdim\lastskip<\bigskipamount\removelastskip\fi\smallskip
            \noindent {\bf Proof.}\quad}

\def\Axiom #1. {\bigbreak\vskip-\parskip\noindent{\bf Axiom #1.}\quad\it}

\def\Satz #1. {\bigbreak\vskip-\parskip\noindent{\bf Satz #1.}\quad\it}

\def\Korollar #1. {\bbr\vskip-\parskip\nin{\bf Korollar #1.} \quad\it}

\def\Bemerkung #1. {\rm\bigbreak\vskip-\parskip\noindent{\bf Bemerkung #1.}
\quad}

\def\Beispiel #1. {\rm\bigbreak\vskip-\parskip\noindent{\bf Beispiel #1.}\quad}
\def\Aufgabe #1. {\rm\bigbreak\vskip-\parskip\noindent{\bf Aufgabe #1.}\quad}

\def\Beweis#1. {\rm\par\ifdim\lastskip<\bigskipamount\removelastskip\fi
           \smallskip\noindent {\bf Beweis.}\quad}

\nopagenumbers

\def\date{\ifcase\month\or January\or February \or March\or April\or May
\or June\or July\or August\or September\or October\or November
\or December\fi\space\number\day, \number\year}

\def\title{Title ??}
\def\author{Author ??}

\def\thanks#1{\footnote*{\eightrm#1}}

\def\rightheadline{\hfil{\eightrm\title}\hfil\tenbf\folio}
\def\leftheadline{\tenbf\folio\hfil{\eightrm\author}\hfil}
\headline={\vbox{\line{\ifodd\pageno\rightheadline\else\leftheadline\fi}}}

\def\firstheadline{}
\def\firstfootline{\cen{\rm\folio}}

\def\seite #1 {\pageno #1
               \headline={\ifnum\pageno=#1 \firstheadline
               \else\ifodd\pageno\rightheadline\else\leftheadline\fi\fi}
               \footline={\ifnum\pageno=#1 \firstfootline\else{}\fi}}

%%%THIS IS THE MACRO LEFTSPACE.TEX %%%TO THD VIA WAFRUPP
\newdimen\dimenone
 \def\checkleftspace#1#2#3#4{%DIESER MACRO STAMMT VON APPELT
 \dimenone=\pagetotal%#1=Skip vorher,#2=Font,#3=Text,#4=Skip nachher  
 \advance\dimenone by -\pageshrink   %testen ob Titel noch mit Gewalt auf Seite 
                                                                          %geht
 \ifdim\dimenone>\pagegoal          %nacha tua nix-- gewoehnliche Outputroutine 
   \else\dimenone=\pagetotal
        \advance\dimenone by \pagestretch
        \ifdim\dimenone<\pagegoal
          \dimenone=\pagetotal
          \advance\dimenone by#1         %addieren Skip vor Ueberschrift (=#1)
          \setbox0=\vbox{#2\parskip=0pt                %#2 ist gewaehlter Font
                     \hyphenpenalty=10000
                     \rightskip=0pt plus 5em
                     \noindent#3 \vskip#4}    %#3=Ueberschrift,#4=skip nachher
        \advance\dimenone by\ht0
        \advance\dimenone by 3\baselineskip   
        \ifdim\dimenone>\pagegoal\vfill\eject\fi
          \else\eject\fi\fi}

%%% OUR HEADLINE MACROS LOOK LIKE THIS USING THIS MACRO

\def\subheadline #1{\nin\bigbreak\vskip-\lastskip
      \checkleftspace{0.7cm}{\bf}{#1}{\medskipamount}
          \indent\vskip0.7cm\centerline{\bf #1}\medskip}

\def\sectionheadline #1{\bigbreak\vskip-\lastskip
      \checkleftspace{1.1cm}{\bf}{#1}{\bigskipamount}
         \vbox{\vskip1.1cm}\cen{\bfone #1}\bsk}

\def\lsectionheadline #1 #2{\bigbreak\vskip-\lastskip
      \checkleftspace{1.1cm}{\bf}{#1}{\bigskipamount}
         \vbox{\vskip1.1cm}\cen{\bfone #1}\msk \cen{\bfone #2}\bsk}

\def\lchapterheadline #1 #2{\bigbreak\vskip-\lastskip\indent\vskip3cm
                       \cen{\bftwo #1} \msk \cen{\bftwo #2} \gsk}
\def\llsectionheadline #1 #2 #3{\bigbreak\vskip-\lastskip\indent\vskip1.8cm
\cen{\bfone #1} \msk \cen{\bfone #2} \msk \cen{\bfone #3} \nobreak\bsk\nobreak}

%\def\[#1 #2\par{\hbox{\vtop{\hsize = 2.5 true cm \nin [#1]\hfill}
%\vtop{\hsize = 12.0 true cm \nin #2\penalty10000\llap.}}
%\vbox{\vskip.3\baselineskip}}  
% parameters for hsize are percentages !!! f.e. 0.2 + 0.8 = 1.0

\newtoks\literat
\def\[#1 #2\par{\literat={#2\unskip.}%
\hbox{\vtop{\hsize=.15\hsize\nin [#1]\hfill}
\vtop{\hsize=.82\hsize\nin\the\literat}}\par
\vskip.3\baselineskip}

\mathchardef\emptyset="001F 
\def\address{Author: \tt$\backslash$def$\backslash$address$\{$??$\}$}

\def\firstpage{\nin
{\obeylines \parindent 0pt }
\vskip2cm
\centerline {\bfone \title}
\gsk
\centerline{\bf\author}

\vskip1.5cm \rm}

\def\addresstwo{}

\def\dlastpage{\par\vbox{\vskip1cm\nin
\line{
\vtop{\hsize=.5\hsize{\parindent=0pt\baselineskip=10pt\nin\address}}
\quad 
\vtop{\hsize=.42\hsize\nin{\parindent=0pt
\baselineskip=10pt\addresstwo}}
\hfill} }}

% END OF LIEMACS.TEX

\def\onto{\to\mskip-14mu\to} 

\pageno=1

\def\title{Unitary spherical highest weight representations}
\def\author{Bernhard Kr\"otz${}^*$
and Karl--Hermann Neeb${}^{**}$}
\footnote{}{\nin ${}^*$Part of the work was supported by the Erwin-Schr\"odinger-Institut, Vienna}
\footnote{}{\nin ${}^{**}$ Part of the work was done on a visit 
supported by the Research Institute of The Ohio State University}

\def\date{June 7, 2000}
%\def\rightheadline{\tenbf\folio\hfil{\bf }\quad\eightrm\date}
% F\"ur endg\"ultige Version l\"oschen
%\def\leftheadline{\tenbf\folio\hfil{\tt spherep.tex}\hfil\eightrm\date}

\def\Box #1 { \msk\par\nin 
\centerline{
\vbox{\offinterlineskip
\hrule
\hbox{\vrule\strut\hskip1ex\hfil{\smc#1}\hfill\hskip1ex}
\hrule}\vrule}\msk }

\def\bs{\backslash} 
\def\str{{\frak {str}}}
\def\Str{\mathop{\rm Str}\nolimits}
\def\le{\mathop{\rm le}\nolimits}
\def\mod{\mathop{\rm mod}\nolimits}
\def\SO{\mathop{{\cal S}'_{\oline\Omega}}\nolimits}
\def\PD{\mathop{{\cal PD}(V,F)}\nolimits}

\def\address
{Bernhard Kr\"otz

The Ohio State University

Department of Mathematics

231 West 18th Avenue

Columbus, OH 43210-1174

USA

kroetz@math.ohio-state.edu 

}

\def\addresstwo
{Karl--Hermann Neeb

Fachbereich Mathematik 

Technische Universit\"at Darmstadt 

Schlo\3gartenstr.\ 7  

D-64289 Darmstadt

Germany

neeb@mathematik.tu-darmstadt.de 

}

\firstpage 

\sectionheadline{Introduction}

\nin Let $G/H$ be a semisimple symmetric space attached to an involution $\tau\: G\to G$. Then 
in order for a unitary 
irreducible representations $(\pi, {\cal H})$ of $G$ to be realized in $L^p(G/H)$, 
$1\leq p\leq \infty$, or more generally in ${\cal D}'(G/H)$, 
it is necessary that $(\pi, {\cal H})$ is $H$-{\it spherical}, i.e., the space of $H$-invariant 
distribution vectors $({\cal H}^{-\infty})^H$ has to be non-zero. 

\par Two classes of representations in the unitary dual $\hat G$ of $G$ are of special interest: The unitary 
principal series and the discrete series. For the unitary principal series induced from a parabolic 
subgroup we have for almost all parameters
a complete description of $({\cal H}^{-\infty})^H$ (cf.\ [Ba88], [\'Ol87] for 
the $\theta\tau$-stable  minimal parabolics and [BrDe92] for the general case). Discrete series on 
$G/H$ 
were constructed in [FJ80] (see also [MaOs84]) and in [Bi90] it was shown that $({\cal H}^{-\infty})^H$ 
is one-dimensional for all discrete series on $G/H$ except for four types of exceptional symmetric spaces. 
Holomorphic discrete series, i.e., unitary highest weight representation of $G$ which 
can be realized in $L^2(G/H)$ were constructed in [\'O\O 88, 91] (this construction is essentially 
different than the one in [FJ80]) and it was shown that 
$({\cal H}^{-\infty})^H$ is always one-dimensional (this also covers half of the exceptional spaces in 
[Bi90]). 

In this paper we address the problem of classifying all $H$-spherical unitary highest weight 
representations for a simply connected hermitian Lie group $G$. We write 
$\g=\h\oplus\q$ for the $\tau$-eigenspace decomposition of $\g\:={\rm Lie} (G)$ and choose  
a $\tau$-stable Cartan decomposition $\g=\k\oplus\p$. Set $K\:=\exp(\k)$. Note that $\z(\k)$ is one-dimensional since 
$G$ is hermitian and that $H$ is connected since 
$G$ is simply connected. 
Now non-trivial $H$-spherical highest weight representations exist if and only if 
$\z(\k)\subeq \q$ which means that $G/H$ is {\it compactly causal} (cf.\ [Hi\'Ol96]). 

\par The unitary highest weight representations $(\pi_\lambda, {\cal H}_\lambda)$ of $G$ 
can be parametrized by the highest weight $\lambda$ of the corresponding finite-dimensional 
unitary representation $(\pi_\lambda^K, F(\lambda))$ of $K$. 
In the {\it scalar case}, i.e., $\dim F(\lambda) = 1$, 
the classification of unitary highest weight representations was accomplished 
in [Wal79] and [VR76] while the general case has been treated in [EHW83] and [Jak83]. 
The  generalized Verma module
$N(\lambda)={\cal U}(\g_\C)\otimes_{{\cal U}(\k_\C+\p^+)}F(\lambda)$ is a highest weight 
module for $\g$ with highest weight $\lambda$. It has a unique irreducible 
quotient $L(\lambda)$ which is $(\g,K)$-isomorphic to the space of $K$-finite 
vectors of $(\pi_\lambda, {\cal H}_\lambda)$. We call $\lambda$ singular if 
the kernel $J(\lambda)$ of the natural map $N(\lambda) \to L(\lambda)$ is non-zero and 
regular otherwise. 

\par If $(\pi_\lambda, {\cal H}_\lambda)$ is spherical, then $L(\lambda)$ trivially is spherical, 
but the converse is not obvious because there is no a priori reason for 
an $\h$-invariant linear functional on $L(\lambda)$ to extend to a continuous functional on 
${\cal H}_\lambda^\infty$. That this is nevertheless the case follows from 
[BaDe88] and [BrDe92] (see also [Kr99a] for the hyperfunction version). 
In the light of this result we only have to deal with the question whether 
$L(\lambda)$ is $\h$-spherical. It is easy to see that a necessary condition is that 
$F(\lambda)$ is $H\cap K$-spherical.  Further elementary 
considerations show that $N(\lambda)$ is $\h$-spherical if and only if $F(\lambda)$ is $H\cap K$-spherical. This already completes the classification for all regular parameters. 
For singular parameters it was first observed in [Kr99a] that the condition  of $F(\lambda)$ being $H\cap K$-spherical 
is not sufficient by showing that the even metaplectic representation of $\Sp(n,\R)$, $n\geq 2$, 
is not $H$-spherical (see also [KN\'O99, Cor.\ II.10]). 

\par If $\lambda$ is singular, then we can attach to $L(\lambda)$ the level ${\rm le}(\lambda)$
which is a natural number ranging from 1 to at most the real rank of $G$. 
Algebraic considerations using the fact that the maximal submodule of $N(\lambda)$ is 
cyclic (cf.\ [DES91]) yield the following half of the classification:  

\msk 

\nin {\bf Theorem II.8.} {\it Suppose that $L(\lambda)$ has odd reduction level. Then 
$(\pi_\lambda, {\cal H}_\lambda)$ is $H$-spherical 
if and only if $F(\lambda)$ is $H\cap K$-spherical. }

\msk    

Besides the algebraic approach to classify $H$-spherical unitary highest weight representations 
one can also attack the problem in an analytic way by using a 
realization of $(\pi_\lambda, {\cal H}_\lambda),$
where the $H$-action is transparent. In the sequel we assume that $G/K$ is of tube type, i.e., 
$G/K$ is biholomorphic to a tube domain $V+i\Omega$ with $V$ an euclidean Jordan algebra and 
$\Omega\subeq V$ the open cone of squares. Then it follows from [HiNe99] that ${\cal H}_\lambda$ 
can be realized as an $L^2$-space of $F(\lambda)$-valued functions on $\oline \Omega$, given 
by 
$$L^2(\oline \Omega, R_\lambda)\:=
\Big\{f\: \oline \Omega\to F(\lambda) \hbox{\rm\ meas.}\: \int_{\oline\Omega}  
\la dR_\lambda(x).f(x), f(x) \ra <\infty\Big\}, $$
where $R_\lambda$ is a certain measure on $\oline \Omega$ 
with values in the cone $\Herm^+(F(\lambda))$ of positive operators on $F(\lambda)$. 
We call $L^2(\oline \Omega, R_\lambda)$ the {\it cone realization}
of $(\pi_\lambda, {\cal H}_\lambda)$. The $H$-action is then given by 
$$(\pi_\lambda(h).f)(x)=\tau_\lambda(\theta(h)).f(\theta(h)^{-1}.x)
\quad \hbox{ for } \quad h\in H, f\in L^2(\oline\Omega,R_\lambda)), $$
where $\tau_\lambda$ denotes the natural $H$-representation on 
$F(\lambda)$ associated to the corresponding representation of $K_\C$ by a 
Cayley transform mapping $H$ into $K_\C$. 
The $K$-action however is completely 
invisible in the cone realization $L^2(\oline\Omega, R_\lambda)$. 
Therefore the following description of the $K$-finite vectors 
is surprisingly explicit and simple: 
\msk
\nin {\bf Theorem IV.5.} {\it The $K$-finite vectors in $L^2(\oline\Omega, R_\lambda)$ 
are given by 
$e^{-\tr}\Pol(V,F(\lambda))\res_{\oline\Omega},$
where $\tr$ denotes the Jordan algebra trace of $V$.}
\msk 
We even obtain a natural identification of $N(\lambda)$ with the space 
$e^{-\tr}\Pol(V,  F(\lambda))$, such that the restriction map to 
$L^2(\oline\Omega, R_\lambda)$ corresponds to the quotient map $N(\lambda) \to L(\lambda)$. 
Using Theorem IV.5  and the results of [HiNe99], we obtain a description of the 
(up to scalar multiples) unique $H$-invariant $F(\lambda)$-valued tempered distribution 
$R_\lambda^H$ supported by $\oline\Omega$. Now $L(\lambda)$ is spherical 
if and only if $R_\lambda^H$ vanishes on the maximal submodule $J(\lambda)$ of $N(\lambda)$, 
realized as $e^{-\tr}\Pol(V,  F(\lambda))$. This requirement leads to the 
necessary condition $\supp(R_\lambda^H) \subeq \supp(R_\lambda)$, but unfortunately 
this condition is not sufficient for $L(\lambda)$ to be spherical. 
Nevertheless, if $\supp(R_\lambda)$ is strictly smaller than $\oline\Omega$, then 
the information on $J(\lambda)$ and $R_\lambda^H$ is explicit enough to obtain the 
following Classification Theorem: 

\msk \nin {\bf Theorem V.11.} {\it Let $G$ be a simply connected hermitian 
Lie group of tube type and $(\pi_\lambda, {\cal H}_\lambda)$ be a unitary 
highest weight representation with $F(\lambda)^{H \cap K} \not=\{0\}$ and 
$\supp(R_\lambda) \not= \oline\Omega$. 
Then $(\pi_\lambda, {\cal H}_\lambda)$ is $H$-spherical if and only if 
${\rm le}(\lambda)$ is odd. }
\msk 

We will see in Section V that among the singular unitary highest weight representations 
with spherical $F(\lambda)$ the condition that $\supp(R_\lambda)$ is smaller than 
$\oline\Omega$ is satisfied in about half of all cases. 
In the scalar case $F(\lambda)$ is always $H\cap K$-spherical and 
$\lambda$ is singular if and only $\supp(R_\lambda) \not=\oline\Omega$. 
Therefore the Classification Theorem covers in particular the scalar case.  
It therefore implies that the most singular scalar type 
highest weight representation which has reduction level $2$ is never spherical. 
This covers the even metaplectic representation of $G = \Sp(n,\R)$. 

In the vector case (i.e. $\dim F(\lambda)>1$) the situation is much more complicated 
if $\supp(R_\lambda)  = \oline\Omega$ and $\lambda$ is singular. Here we 
do not know of any example with even reduction level, where $L(\lambda)$ is spherical. 

\par We conclude the paper with an appendix where we give a simple algebraic 
proof of the fact that parameters $\lambda$ which belong to the relative holomorphic 
discrete series representations of $G/H$ are always regular. 
This means that these representations do not provide any further information on the question 
whether a singular unitary highest weight representation $L(\lambda)$ with spherical lowest 
$K$-type $F(\lambda)$ is spherical or not. 

\sectionheadline{I. Generalities on highest weight modules} 

\nin In this section we recall some concepts and results 
related to spherical highest weight modules. 

\subheadline{Hermitian Lie algebras}
 
Let $\g$ be a simple real Lie algebra and $\theta$ a Cartan involution of $\g$ 
with Cartan decomposition $\g=\k\oplus\p$. 
We call $\g$ {\it hermitian} if  $\z(\k)\neq \{0\}$. 
We now collect some basic facts concerning hermitian Lie algebras (cf.\ [Hel78, 
Ch.\ VIII]). 
The center of $\k$ is one-dimensional, i.e., $\z(\k)=\R i Z_0$ for some 
$0\neq Z_0\in i \z(\k)$. We can normalize $Z_0$ such that $\Spec (Z_0)=\{ 1,0,-1\}$. Further, 
every Cartan subalgebra $\t$ of $\k$ is a Cartan subalgebra of $\g$ and 
$\z(\k)\subeq \t$. Let $\g_\C$ be the complexification of $\g$ and $\hat\Delta$ the 
root system of $\g_\C$ with respect to $\t_\C$. 

\par A root $\hat\alpha\in \hat\Delta$ is called {\it compact} if $\hat\alpha(Z_0)=0$   
and {\it non-compact} otherwise. We denote by $\hat\Delta_k$, resp. 
$\hat\Delta_n$, 
the set  of compact, resp. non-compact roots. We fix a positive system 
$\hat\Delta^+\subeq \hat\Delta$ such that 
$$\hat\Delta_n^+\:=\hat\Delta_n\cap \hat\Delta^+
=\{ \hat\alpha\in \hat\Delta_n\: \hat\alpha(Z_0)=1\}.$$
We set $\p^\pm\:=\{X\in\g_\C\: [Z_0,X]=\pm X\}$ and note that 
$$\g_\C=\p^+\oplus\k_\C\oplus\p^-.$$
As $\Spec(Z_0)=\{1,0,-1\}$ and $Z_0\in \z(\k_\C)$, it follows that 
$[\k_\C, \p^\pm]\subeq \p^\pm$, $[\p^+, \p^-]\subeq \k_\C$, $[\p^+, \p^+]=\{0\}$ and 
$[\p^-, \p^-]=\{0\}$.

\subheadline{Generalized Verma  modules}

In this subsection we collect some basic facts concerning highest weight modules 
from an abstract algebraic point of view. References for the fact used in this subsection 
are [EHW83] or [Ne99, Ch.\ IX].  

\par Let $X\mapsto\oline X$ denote the conjugation in $\g_\C$ with respect to the 
real form $\g$. Then the map $X\mapsto X^*\:=-\oline X$ extends to an involutive 
antilinear antiautomorphism $D\mapsto D^*$ of ${\cal U}(\g_\C)$.
A hermitian form $\la\cdot , \cdot\ra$ on a $\g_\C$-module $V$ is called 
{\it contravariant} if 
$$(\forall X\in \g_\C)(\forall v, w\in V)\quad 
\la X.v, w\ra =\la v, X^*.w\ra. $$ 

Let $\lambda\in i\t^*$ be dominant integral with respect to 
$\hat\Delta_k^+$ and write 
$F(\lambda)$ for the corresponding finite-dimensional unitary irreducible $\k_\C$-module. Let 
$\q:=\p^+\rtimes \k_\C$
and turn $F(\lambda)$ into a $\q$-module by letting $\p^+$ act trivially. 
We define the {\it generalized Verma module} associated to $\lambda$ by 
$$N(\lambda)={\cal U}(\g_\C)\otimes_{{\cal U}(\q)} F(\lambda).$$
Then $N(\lambda)$ is a highest weight module  
with respect to $\hat\Delta^+$ and highest weight $\lambda$ which contains $F(\lambda)$ 
as a $\k$-submodule which is called the {\it lowest $\k$-type}. 
Therefore it has a unique maximal submodule $J(\lambda)$ and hence a unique 
simple quotient $L(\lambda) := N(\lambda)/J(\lambda)$. 
The positive definite form on $F(\lambda)$ has a unique extension to 
$N(\lambda)$ called the {\it Shapovalov form}. 
The radical of this form is the 
maximal submodule $J(\lambda)$. In particular, the Shapovalov form factors to a 
contravariant form on 
$L(\lambda)$ which we also denote by $\la \cdot, \cdot\ra_\lambda$.  
We call $L(\lambda)$ 
{\it unitarizable} if $\la\cdot, \cdot\ra_\lambda$ is positive definite on 
$L(\lambda)$.   

\par  The canonical mapping 
${\cal U}(\p^-)\otimes F(\lambda)\to N(\lambda)$ 
gives rise to an isomorphism 
$${\cal U}(\p^-)\otimes F(\lambda)\cong N(\lambda)$$
of $\k_\C\ltimes\p^-$-modules, where 
$\k_\C\ltimes \p^-$ acts on ${\cal U}(\p^-)\otimes F(\lambda)$ 
via
$$\eqalign{&(\forall X\in \k_\C)\qquad  X.(p\otimes v)
\:= X.p\otimes v +p\otimes X.v\cr 
&(\forall Y\in \p^-)\qquad  Y.(p\otimes v)\:=Yp\otimes v \cr}\leqno(1.1) $$
for all $p\in {\cal U}(\p^-)$ and $v\in F(\lambda)$.

\subheadline{Symmetric structures}

We now endow $\g$ with a symmetric structure, i.e., with an 
involutive automorphism $\tau\: \g\to\g$. The pair 
$(\g,\tau)$ is  called a {\it symmetric Lie algebra}. We write 
$$\h\:=\{X\in\g\: \tau(X)=X\}\quad\hbox{and}\quad 
\q\:=\{ X\in\g\: \tau(X)=-X\}$$
for the $\tau$-eigenspaces and note that $\g=\h\oplus\q$. 
We assume that  $\tau$ commutes with the Cartan involution $\theta$.

\Definition I.1. A hermitian symmetric Lie algebra $(\g,\tau)$
is called {\it compactly causal}, if $\z(\k)\subeq \q$.\qed 

\Remark I.2. For 
later reference we collect here the basic structural facts concerning 
compactly causal symmetric Lie algebras. 
Let $(\g,\tau)$ be compactly causal. We denote the complex linear 
extensions of $\tau$ and $\theta$ to $\g_\C$ by the same symbols. 
The {\it c-dual} $(\g^c,\tau^c)$ of $(\g,\tau)$ is defined by 
$\g^c\:=\h+i\q$ and $\tau^c\:=\tau\res_{\g^c}$. Then 
$\theta^c\:=\theta\tau\res_{\g^c}$ defines a Cartan involution 
on $\g^c$ and we write $\g^c=\k^c\oplus\p^c$ for the corresponding 
Cartan decomposition. 
\par As $(\g,\tau)$ is compactly causal, we find a maximal 
abelian subspace $\a\subeq (i\q)\cap\p^c$ which is already 
maximal abelian in $\p^c$ (cf.\ [Hi\'Ol96, Prop.\ 3.1.11]). We write 
$\Delta=\Delta(\g^c,\a)$ for the restricted root system of 
$\g^c$ with respect to $\a$ and write 
$$\g^c=\a\oplus\z_\h(\a)\oplus\bigoplus_{\alpha\in\Delta} (\g^c)^\alpha$$
for the corresponding root space decomposition.  
We now choose
the compact Cartan subalgebra 
$\t$ of $\g$ such that $\t=(\t\cap\h)\oplus i\a$. Since $\z(\k)\subeq \q$ we can 
choose $\hat\Delta^+$ such that the prescription 
$\Delta^+\:=\hat\Delta^+\res_\a\bs\{0\}$ defines a positive system $\Delta^+$ of $\Delta$. 
Similarly we define {\it compact } and {\it non-compact roots} 
$\Delta_k$ and $\Delta_n$. Finally we define subalgebras 
$$\n^+\:=\bigoplus_{\alpha\in \Delta^+} (\g^c)^\alpha
\qquad\hbox{and}\qquad \n^-\:=\bigoplus_{\alpha\in -\Delta^+} (\g^c)^\alpha
.$$ 
Note that $\g^c=\h\oplus\a\oplus\n^+$ and $\n^\pm=\p^\pm\cap\g^c$.\qed

\subheadline{Spherical highest weight modules}

\Definition I.3. (a) If $\h$ is a Lie algebra and 
$V$ a complex $\a$-module, then we write 
$V^\h$ for the submodule of $V$ which is annihilated by $\h$. 
We denote by $V^\sharp$ the algebraic antidual of $V$, i.e., the space of 
antilinear functionals $V \to \C$.

\par\nin (b) If $(\g,\tau)$ is a hermitian symmetric 
Lie algebra and $L(\lambda)$ an irreducible highest weight module, then 
we call $L(\lambda)$ {\it spherical} if there exists 
$0\neq \nu\in L(\lambda)^\sharp$ which is annihilated by $\h$, i.e., 
$(L(\lambda)^\sharp)^\h\neq\{0\}$. \qed    

\Proposition I.4. If the 
hermitian symmetric Lie algebra $(\g,\tau)$ has a 
non-trivial spherical unitarizable highest weight 
module, then $(\g, \tau)$ is compactly causal. 

\Proof.  [Kr99b, Lemma B.1].\qed 

\Proposition I.5. Let  $V$ be a highest weight module 
of the compactly causal symmetric Lie algebra $(\g, \tau)$ 
with highest weight $\lambda$ and let $v_\lambda$ denote 
a highest weight vector. Then 
$${\cal U}(\h_\C).v_\lambda= V.$$
If, in addition, $V$ is spherical, then the  lowest $\k_\C$-type 
$F(\lambda) := \{ v \in V \: \p^+.v = \{0\}\}$ 
of $V$ is spherical for $(\k,\tau)$, i.e., 
$F(\lambda)^{\h\cap\k}\neq \{0\}$. 

\Proof. Since $(\g,\tau)$ is compactly causal, 
the decomposition $\g^c=\h\oplus\a\oplus\n$ implies that 
$\g_\C=\h_\C\oplus\a_\C\oplus\n_\C^+$. Therefore 
${\cal U}(\g_\C)={\cal U}(\h_\C){\cal U}(\a_\C){\cal U}(\n_\C^+)$
by the Poincar\'e-Birkhoff-Witt Theorem, and therefore 
$$ V={\cal U}(\g_\C).v_\lambda=
{\cal U}(\h_\C){\cal U}(\a_\C){\cal U}(\n_\C^+).v_\lambda=
{\cal U}(\h_\C).v_\lambda. $$

Now we assume that $V$ is spherical and consider 
$\nu \in (V^\sharp)^\h \setminus \{0\}$. 
Then the first part implies that $\la\nu,v_\lambda\ra\neq0$. Thus
$\nu\res_{F(\lambda)}$ defines a non-zero $\h\cap\k$-fixed element.
\qed 

The following lemma provides some refined information on the $\h$-module 
structure of $N(\lambda)$. 

\Lemma I.6. If $(\g,\tau)$ is compactly causal, then 
we have an isomorphism of ${\cal U}(\h_\C)$-modules: 
$${\cal U}(\h_\C)\otimes_{{\cal U}(\h_\C\cap\k_\C)} F(\lambda)
\to N(\lambda)= {\cal U}(\g_\C)\otimes_{{\cal U}(\k_\C\oplus \p^+)} 
F(\lambda), \ \ u\otimes v\mapsto u\otimes v.$$

\Proof. [HiKr98, Lemma 3.1.1].\qed 

\Proposition I.7. If $(\g,\tau)$ is 
compactly causal, then the restriction map 
$$ (N(\lambda)^\sharp)^\h \to (F(\lambda)^\sharp)^{\h\cap\k}, 
\quad \nu \mapsto \nu\res_{F(\lambda)} $$
is a bijection. 
In particular, $\dim (N(\lambda)^\sharp)^\h \leq 1$ and 
$N(\lambda)$ is spherical if and only if $F(\lambda)$ is. 
If $F(\lambda)$ is spherical, then $L(\lambda)$ is spherical if and only 
a non-zero element $\nu \in (N(\lambda)^\sharp)^\h$ 
vanishes on the kernel $J(\lambda)$ of the map 
$N(\lambda) \onto L(\lambda)$. 

\Proof. (cf.\ [HiKr98, Lemma 3.1.2]) This is a direct consequence of 
Lemma I.6.
\qed

\sectionheadline{II. The classification for odd reduction levels}

\nin We are interested in a description of those unitary highest weight modules
$L(\lambda)$ which are spherical. In Proposition I.7 we have seen that 
it is necessary that $F(\lambda)$ is spherical, and that this is sufficient 
provided that $N(\lambda) \cong L(\lambda)$. Therefore the question becomes 
interesting for the situations where the kernel $J(\lambda)$ of the 
quotient map $N(\lambda) \onto L(\lambda)$ is non-trivial. In this 
case we will assign to $L(\lambda)$ a natural number $\le(\lambda)$, called 
its reduction level. 
The main result of this section is that if 
$\le(\lambda)$ is odd and $F(\lambda)$ is spherical, then $L(\lambda)$ is 
spherical. This is done by analyzing the $\t_\C$-Fourier series of 
the spherical vector $\nu$ of $L(\lambda)$. In the 
following sections we will take a closer look at the cases 
with even reduction level. 

\subheadline{The classification of unitarizable highest weight modules}  

Let $\g$ be a hermitian Lie algebra, 
$\t \subeq \k$ a Cartan subalgebra, and identify $\z(\k)^*$ with the subspace 
$(\t \cap [\k,\k])^\bot \subeq \t^*$. Let 
$\zeta\in i\z(\k)^* \subeq i\t^*$ be normalized by 
$\zeta(\check \gamma)=1$ for $\gamma \in \hat\Delta_n^+$ a long root. 
Fix $\lambda^0\in i\t^*$  dominant integral with respect 
to $\hat \Delta_k^+$. For $u\in \R$ we set 
$$\lambda^u=\lambda^0-u\zeta.$$
Further we define
$$l(\lambda^0)\:=\{ u\in \R \: L(\lambda^u)\  \hbox{is unitarizable}\}.$$ 

\Theorem  II.1. {\rm (Enright--Howe--Wallach, Jakobsen)}
Let $\lambda^0\in i\t_0^*$ be dominant integral with 
respect to $\hat\Delta_k^+$. 
\item{(i)} There exists a real number $u_m$ such that 
$$l(\lambda^0)=\{ u_0, \ldots, u_m\} \dot\cup ]u_m,\infty[ $$
where $u_j = u_0 + j{d\over 2}$ for $j =0,\ldots, m$ and 
$d$ is the multiplicity of the restricted roots of $\g$ which are the second longest. 
\item{(ii)} For $u\in l(\lambda^0)$ we have 
$L(\lambda^u)\cong N(\lambda^u)$ if and only if $u> u_m$. 

\Proof. [EHW83, Th.\ 2.4]. \qed 

The eigenspace decomposition of the 
symmetric algebra ${\cal U}(\p^-)$ with respect to the 
action of $Z_0$ by derivations yield the grading by homogeneous elements  
${\cal U}(\p^-)=\bigoplus_{n=0}^\infty {\cal U}(\p^-)^n$. For each 
$n\in \N_0$ we set 
$N(\lambda)^n\:={\cal U}(\p^-)^n\otimes F(\lambda).$
Then $N(\lambda)=\bigoplus_{n=0}^\infty N(\lambda)^n$ as $\k_\C$-modules, 
and (1.1) shows that $Z_0$ acts on 
$N(\lambda)^n$ by multiplication with $\lambda(Z_0) - n$.

\Definition II.2.  Suppose that $\lambda\in i\t^*$ is dominant integral with respect to 
$\hat\Delta_k^+$ and let $J(\lambda)$ be the maximal submodule of $N(\lambda)$. 
We call $\lambda$ {\it regular} if $J(\lambda)=\{0\}$ and {\it singular} otherwise. 
Let $J(\lambda)=\bigoplus_{n=0}^\infty J(\lambda)^n$ be the natural grading 
induced from the one on $N(\lambda)$. If $\lambda$ is singular, then we call the number  
$${\rm le}(\lambda)\:=\min\{n\in\N_0\: J(\lambda)^n\neq \{0\}\}$$
the {\it level of reduction } of $L(\lambda)$.
\qed 

\Remark II.3.  From the classification of unitary 
highest weight modules, we know that the points $u_0,\ldots, u_m$ from Theorem II.1 are 
related to reduction levels by 
$$ \le(\lambda^{u_j}) = 1 + j $$
(cf.\ [EJ90]). 
So we see that approximately half of the singular points $u_0,\ldots, u_m$ correspond to 
even and half to odd reduction levels. 
\qed

\Theorem II.4. {\rm (Davidson--Enright--Stanke)} 
Suppose that $N(\lambda)$ is reducible and 
$L(\lambda)$ is unitarizable. Then the maximal submodule $J(\lambda)$ is a highest weight 
module, i.e., there exist an element $\lambda_J\in i\t^*$ which is dominant integral 
with respect to $\hat\Delta_k^+$ such that $J(\lambda) \cong N(\lambda_J)$.

\Proof. This follows from [DES91, Th.\ 3.1]. Note that, since
$N(\lambda)$ is a free ${\cal U}(\p^-)$-module, the maximal submodule
$J(\lambda)$ is isomorphic to $N(\lambda_J)$ if $\lambda_J$ is its
highest weight. 
\qed   

\subheadline{Some algebraic criteria}

Let 
$\sigma\:=\exp(\pi i Z_0)\in G$. Then the Cartan involution $\theta$ on $G$ is given by 
$\theta(g)=\sigma g\sigma^{-1}$ for $g\in G$.

\Lemma II.5. {\rm (Parity Lemma)} Let $N(\lambda)$ be spherical and 
$\nu\in (N(\lambda)^\sharp)^\h$ be a
non-zero spherical vector. For $n\in\N_0$ we
set $\nu^n\:=\nu\res_{N(\lambda)^n}$. Then 
$\nu^n=0$ for $n$ odd. 

\Proof. Since the action of $\k$ on $N(\lambda)$ is locally finite, it integrates a 
representation of the simply connected group $K \subeq G$. Hence, in particular, 
$\sigma$ acts on $N(\lambda)$ in such a way that 
$\sigma(X.v) = \theta(X).(\sigma.v)$ for $X \in \g$, $v \in N(\lambda)$. 
Since $\h$ is $\theta$-invariant, $\nu \circ \sigma = c \nu$ for some 
$c\in \C^\times$ follows from 
$(N(\lambda)^\sharp)^{\h \cap \k} = \C \nu$ (Proposition I.7). 
From $\nu=\sum_{n=0}^\infty \nu^n$ and (1.1) we thus get 
$$\sigma.\nu=e^{\pi i \lambda(Z_0)} \sum_{n=0}^\infty (-1)^n \nu^n.$$
Therefore $c=e^{\pi i \lambda(Z_0)}$ and $\nu^n=0$ for $n$ odd.
\qed 

{}From now on we suppose that $(\g,\tau)$ is compactly causal, that 
$L(\lambda)$ is unitarizable, 
and that $F(\lambda)$ is spherical, so that 
$N(\lambda)$ is spherical by Proposition I.7. 
Let $0\not= \nu \in (N(\lambda)^\sharp)^\h$. 

\Lemma II.6. If $J(\lambda)$ is not spherical, then $L(\lambda)$ is spherical.
 
\Proof. If $J(\lambda)$ is not spherical, 
then $\nu\res_{J(\lambda)}=0$. 
Thus $\nu$ factors to a non-zero $\h$-fixed antilinear functional 
on $L(\lambda)$, showing that $L(\lambda)$ is spherical. \qed

\Proposition II.7. Let $0\neq\nu\in (N(\lambda)^\sharp)^\h$ and 
$F(\lambda_J)$ the highest $\k$-type of $J(\lambda)$. Then the following 
assertions are equivalent: 
\item{(1)} $L(\lambda)$ is spherical. 
\item{(2)} $\nu\res_{F(\lambda_J)}=0$. 

\Proof. Note that $L(\lambda)$ is spherical if and only if 
$\nu\res_{J(\lambda)}=0$. We 
recall from Theorem II.4 that $J(\lambda)$ is a highest weight module 
with highest weight $\lambda_J$. Thus Proposition I.5 implies 
that $\nu\res_{J(\lambda)}=0$ if and only if $\nu\res_{F(\lambda_J)}=0$. 
\qed 

\Theorem II.8. {\rm (Classification for odd reduction levels)}
If $L(\lambda) \not\cong N(\lambda)$ 
is unitarizable and $\le(\lambda)$ is odd, then 
$L(\lambda)$ is spherical if and only if $F(\lambda)$ is spherical. 

\Proof. If $L(\lambda)$ is spherical, then $F(\lambda)$ is spherical by 
Proposition I.5. 
\par Assume, conversely, that $F(\lambda)$ is spherical. Then 
$N(\lambda)$ is spherical (Proposition I.7), so that there exists 
a non-zero element  $\nu\in (N(\lambda)^\sharp)^\h$. In view of 
Proposition II.7,  we have to show that $\nu\res_{F(\lambda_J)}=0$. 
Since ${\rm le}(\lambda)$ is odd, this follows from 
Lemma II.5.\qed

\Example II.9. In the following we give an example where both $L(\lambda)$ and $J(\lambda)$
are spherical. 
\par Let $(\g,\h)=(\sp(n,\R),\gl(n,\R))$. We choose $\lambda^0=0$ and set 
$\lambda\:=\lambda^{u_{n-1}}$ for the first reduction point which has reduction 
level $n$. Then it follows from [EJ90] that the highest weight $\lambda_J$ of the 
maximal submodule $J(\lambda)$ is contained in $i\z(\k)^*$. In particular, 
$F(\lambda_J)$ is $\h\cap\k$-spherical. Since $\lambda$ is a first reduction point, 
we have $J(\lambda)=N(\lambda_J)$. Hence $J(\lambda)$ is spherical by Proposition I.7. 
Further if $n$ is odd, then Theorem II.8 implies that $L(\lambda)$ is 
spherical too. More details on the composition series of the modules $N(\lambda)$ in this 
case can be found in [FaKo90]. 
\qed

\sectionheadline{III. Realization of representations in spaces of distributions} 

\nin In this section we realize simple highest weight modules of 
the important special case of Cayley type 
compactly causal symmetric Lie algebras in spaces of vector-valued distributions 
supported by a closed convex cone in a euclidean Jordan algebra. 
For the constructions in this section unitarizability will play no role at all. 
In the next section we will see how this picture can be used to get refined information 
on unitary highest weight modules because in this case the corresponding realization 
by distributions can also be viewed as a subspace of a suitable $L^2$-space with respect to 
an operator-valued measure. 

\subheadline{Algebraic preliminaries}

\Definition III.1. A compactly causal symmetric Lie algebra $(\g,\tau)$ is said to be of 
{\it Cayley type} if $\z(\h)\neq \{0\}$. \qed 

\Remark III.2. (a) If $(\g,\tau)$ is of Cayley type, then $\z(\h)\subeq \h\cap\p$ 
and $\dim \z(\h)=1$ (cf.\ [Hi\'Ol96, Th.\ 1.3.11]).  
\par\nin (b)  Cayley type symmetric spaces are classified (cf.\ [Hi\'Ol96, Th.\ 3.2.8]). 
Up to isomorphy, the corresponding pairs $(\g,\h)$ are given by 
$$(\su(n, n), \sL(n, \C)\oplus\R), \ (\so^*(4n), \su^*(2n)\oplus\R),\ (\sp(n,\R), \gl(n,\R)),$$
$$(\so(2,n), \so(1,n-1)\oplus\R),\ (\e_{7(-25)}, \e_{6(-26)}\oplus\R).$$

\par\nin (b) If $(\g,\tau)$ is of Cayley type, then the Lie algebra 
$\g$ is hermitian and of {\it tube type} (cf.\ [Hi\'Ol96], [KoWo65]), i.e., the
hermitian symmetric space $G/K$ attached to $\g$ is biholomorphically 
equivalent to a tube domain $T_\Omega=V +i \Omega$ over the convex open cone $\Omega$ 
of invertible squares in a finite-dimensional simple euclidean Jordan algebra $V$. 
\qed 

In the following $V$ denotes a simple euclidean Jordan algebra with unit element $e$ and 
$\Omega \subeq V$ the open cone of invertible squares. 
Having Remark III.2(b) in mind, from now on $G$ stands for the simply connected covering group 
of $\Aut(T_\Omega)_0$ and $K = \{ g \in G \: g.ie = ie \}$ is the analytic subgroup 
corresponding to $\k$. The corresponding Cartan involution $\theta$ satisfies 
$\eta(g.\eta^{-1}(z)) = \theta(g).z$ for $z \in T_\Omega$, 
where $\eta(z)=-z^{-1}$. As a subgroup of $G$, the group $H$ consists of all those 
elements of $G$ acting by linear maps on $T_\Omega$. This implies in particular 
that $H$ preserves $V$ and also $\Omega$. Geometrically the involution $\tau$ is 
determined by $\zeta(g.\zeta^{-1}(z)) = \tau(g).z$ for $z \in T_\Omega$, where 
$\zeta(z) = - \oline z$. Since $G$ is simply connected, the subgroup 
$G^\tau$ is connected, hence coincides with $H$. 
Let $K_\C$ denote the universal complexification of $K$. Since 
$K_\C$ has a polar decomposition $K_\C = K\exp(i\k)$, this group is also simply connected. 

\Definition III.3. (a) We follow the notation of [FaKo94]. We 
define $L(z)\in \End(V_\C)$, $z\in V_\C$, by $L(z).w=zw$ and write 
$$P\: V_\C\to \End(V_\C),\ \  z\mapsto P(z)=2L(z)^2-L(z^2)$$ 
for the 
the quadratic representation of $V_\C$. 
Let $\{c_1,\ldots, c_r\}$ be a 
Jordan frame of $V$ (cf.\ [FaKo94, p.\ 44]) 
and set $e_k\:=c_1+\ldots+c_k$ for $0\leq k\leq r$ with the 
conventions $e_0=0$ and $e = e_r$. Then 
$\Omega=H.e_r$  and 
$$\oline \Omega={\cal O}_0 \dot\cup {\cal O}_1\dot\cup\ldots
\dot \cup {\cal O}_r$$
with ${\cal O}_j= H.e_j$ and $\oline {{\cal O}_j}=\cup_{i<j} 
{\cal O}_i$ is the $H$-orbit decomposition of $\oline\Omega$ (cf.\ [FaKo94, Prop.\ IV.3.1]). 
The number $r$ is called the {\it rank of the Jordan algebra $V$}. 

Note that $\h\cap\k \cong \der(V)$ 
and $\h\cap\p \cong \{L(x)\: x\in V\}$ (cf.\ [FaKo94, Th.\ III.5.1]). 
Then $\b :=\oplus_{j=1}^r \R L(c_j)$ is a maximal abelian subspace of $\h \cap \p$. 
Define $\eps_j\in \b^*$ by $\eps_j(L(c_k))=\delta_{jk}$ and set $\eps\:=\sum_{j=1}^r \eps_j$. 
Then the restricted root 
system of $\h$ is given by 
$$ \Sigma = \{ {\textstyle{1\over 2}}(\eps_j - \eps_i) \: i \not= j, 1 \leq i,j \leq r\}.$$
We consider the positive system 
$\Sigma^+ = \{ {\textstyle{1\over 2}}(\eps_j - \eps_i) \: i < j\}.$
The coroot corresponding to ${1\over 2}(\eps_j - \eps_i)$ is 
$2(L(c_j) - L(c_i))$. Therefore a linear functional 
$\lambda = \sum_{j=1}^r {m_j\over 2} \eps_j$ is dominant integral, i.e., 
can occur as a restricted highest weight of a finite-dimensional representation 
of $\h$ if and only if 
$$ m_j - m_i \in \N_0, i < j, $$
which implies in particular that 
$m_1 \leq m_2 \leq \ldots \leq m_r.$
For real highest weight representations which are the same as spherical representations of $H$, 
one has the stronger condition that $m_j - m_i \in 2\N_0$ for $i < j$ 
(cf.\ [HiNe99, App.\ A] and [Hel84, Th.\ 5.1]). 

\par\nin (b) Let 
$$ {\cal D} =\{ z\in V_\C\: I-P(z)P(\oline z) \> 0\} \subeq\p^+ \cong V_\C $$ 
be the Harish Chandra realization of $G/K$ as a bounded symmetric domain 
(cf.\ [FaKo94, Prop.\ X.4.2, Th.\ X.4.3]). 
Then the Cayley transform 
$$c\: {\cal D}\to T_\Omega, \quad z\mapsto -i(z+e)(z-e)^{-1} $$
is biholomorphic. 

The universal complexification $G_\C$ of $G$ acts by meromorphic mappings on 
$V_\C$ and $\g_\C$ can be viewed as a Lie algebra of holomorphic vector fields on $V_\C$.
In this sense the Cayley transform can be considered as obtained by an element of 
$G_\C$, which implies that it induces an isomorphism $c_\g \in \Aut(\g_\C)$ mapping 
the Lie algebra $\tilde \g$ of $\Aut({\cal D})$ onto $\g$. Note that $c_\g=e^{i{\pi\over 2} \ad T}$ for some 
$\h\cap\k$-fixed element $T\in \q\cap\p$ with $\Spec(\ad T)=\{-1,0,1\}$ (cf.\ [Hi\'Ol96]). 
The subalgebra 
$\tilde \k \subeq \tilde \g$ consisting of all vector fields vanishing in $0$ is maximal 
compact, and furthermore we have 
$$ \tilde \k =  (\h \cap \k) \oplus i(\h \cap \p) $$
(cf.\ [FaKo94, Prop.\ X.3.1, Th.\ X.5.3]). 
It follows in particular that $\tilde \k_\C = \h_\C$ and that 
$\tilde \k$ is a compact real form of $\h_\C$. Therefore 
complex $\h$-modules are equivalent to complex $\tilde\k$-modules. 
In addition, we have 
$$ \k =  c_\g(\tilde \k) = (\h \cap \k) \oplus i c_\g(\h \cap \p) $$
and $ic_\g(\h \cap \p)=\q\cap\k$. We choose a compactly embedded Cartan subalgebra $\t \subeq \g$ containing 
$i c_\g(\b)$, which is then automatically adapted to the above decomposition of $\k$. 
We choose a positive system $\hat\Delta_k^+$ of the corresponding root system 
$\hat\Delta_k = \Delta(\k_\C,\t_\C)$ such that 
$\hat\Delta_k^+ \circ c_\g\res_\b \subeq \Sigma^+ \cup \{0\}$. 

If $(\pi_\lambda^\k, F(\lambda))$ is a highest weight representation of $\k$, resp., $\k_\C$ 
we obtain a representation of $\h_\C = \tilde \k_\C$ by 
$\tau_\lambda^\h := \pi_\lambda^\k \circ c_\g$. We call 
$\lambda \circ c_\g\res_\b$ the {\it restricted highest weight} of this representation 
(cf.\ (a)).

On the group level the inclusion 
$c_\g^{-1}\: \h\to\k_\C$ yields the identity on 
$H\cap K$. Therefore the Cartan decomposition $H = (H \cap K)\exp(\h\cap \p)$ 
leads to a homomorphism $H \to K_\C$ which, in view of the simple connectedness of $K_\C$, 
is the universal complexification of $H$, i.e., $H_\C \cong K_\C$. 
\qed

\subheadline{Polynomial differential operators and distributions} 

For a finite-dimensional complex Hilbert vector space $F$, let 
$\SO(F)$ denote the space of $F$-valued tempered distributions 
on $V$ supported by the closed convex cone $\oline\Omega$. 
For $z \in V_\C$ let $e_{z} \: V \to \C$ be defined by 
$e_z(x) := e^{(z,x)}$, where $(\cdot, \cdot) \: V_\C \times V_\C \to \C$
is the complex bilinear extension of the scalar product on $V$. 
According to Lemma B.1 in [HiNe99], for each $z \in T_\Omega$ and 
$D \in \SO(F)$ we have $e_{iz} D \in \SO(F)$ and 
$e_{iz}\res_{\oline\Omega}$ extends to a Schwartz function on $V$, so that the 
Fourier transform ${\cal F}(D)$ extends uniquely to a holomorphic 
function 
$${\cal F}(D) \: T_\Omega \to F, \quad 
z \mapsto  D(e_{iz}), $$
(cf.\ [Sch66, Ch.\ VIII, Prop.\ 6]). 
We thus obtain an injective map 
$$ {\cal F} \: \SO(F) \to \Hol(T_\Omega,F). $$
For a characterization of the image we refer to [Sch66, Ch.\ VIII]. Here 
we will only need that this map is injective. 

In the following we will consider the space
${\cal S}'(V,F)={\cal S}'(V)\otimes V$ of $F$-valued distributions as the 
topological antidual space of ${\cal S}(V,F)$, the $F$-valued Schwartz space, 
where the pairing is given by 
$$ \la D \otimes v, \phi \otimes w \ra := \oline{D(\phi)} \la v, w \ra. $$

Now we turn to the action of the algebra of polynomial differential 
operators on these spaces. To fix the notation, for 
$v \in V$ we define on $C^\infty(V,F)$ the operators 
$$ (\partial_v.\phi)(x) = d\phi(x)(v) \quad \hbox{ and } \quad 
(m_v.\phi)(x) = (v,x) \phi(x). $$
It is also clear that $\End(F)$ acts naturally on $C^\infty(V,F)$ by 
$A.\phi := A \circ \phi$. We write 
${\cal PD}(V,F)$ for the algebra generated by the operators 
$\partial_v$, $m_v$ and $\End(F)$. The elements of this 
algebra are sums of operators of the type 
$P(x) \partial_{v_1} \cdots \partial_{v_n},$
where $P \: V \to \End(F)$ is a polynomial function. 
We define an antilinear involution $*$ 
on ${\cal PD}(V,F)$ by prescribing on the generators that 
$$ m_v^* := m_v, \quad \partial_v^* := - \partial_v \quad \hbox{ and } \quad 
A^* := A^*, $$
where $A^*$ for $A \in \End(F)$ denotes the operator adjoint. Then we 
obtain a natural action of $\PD$ on $\SO(F)$ by 
$$ \la P.D, \phi \ra := \la D, P^*.\phi \ra. $$
Here we use that the elements of $\PD$ restrict to continuous operators 
on the Schwartz space ${\cal S}(V,F)$. In particular we have 
$$ \la \partial_v.D,\phi\ra = \la D,-\partial_v.\phi\ra 
\quad \hbox{ and } \quad 
\la m_v.D,\phi\ra = \la D,m_v. \phi\ra. $$

It is also clear how the algebra $\PD$ acts on the space 
$\Hol(T_\Omega,F)$. In this sense we obtain directly from 
the definitions: 
$$ {\cal F}(\partial_v.D)=  -i m_{v}.{\cal F}(D), \quad 
{\cal F}(m_v.D) =  -i \partial_v.{\cal F}(D) \quad \hbox{ and } \quad 
{\cal F}(A.D) =  A.{\cal F}(D) \leqno(3.1) $$
for $v \in V$ and $A \in \End(F)$. 
The preceding relations describe 
in which way the natural action of the algebra 
${\cal PD}(V,F)$ on $\SO(F)$ is intertwined with its natural 
action on the space $\Hol(T_\Omega,V)$. The Fourier transform 
defines an algebra automorphism 
${\cal F} \: \PD \to \PD$
which satisfies for all $P \in {\cal PD}(V,F)$ and $D\in {\cal S}'(V,F)$ 
$${\cal F}(P.D) = {\cal F}(P) {\cal F}(D). \leqno(3.2) $$
Now (3.1) can be written as 
$$ {\cal F}(\partial_v) = -i m_{v}, \quad 
{\cal F}(m_v) = -i \partial_v \quad \hbox{ and } \quad 
{\cal F}(A) = A. $$

\subheadline{Lie algebra actions} 

Recall that we consider $G$ as the universal covering group of 
$\Aut(T_\Omega)_0$. We write $\Str(V_\C)$ for the structure group of
the complex Jordan algebra $V_\C$ (cf.\ [FaKo94,p.\ 147]). Then the 
natural homomorphism $H \to \Str(V_\C)$, coming 
from the inclusion $H \subeq G$, 
induces a holomorphic covering $H_\C\to \Str(V_\C)_0$, where 
$H_\C$ denotes the universal Lie group complexification of $H$.  
The group $G$ naturally acts on $T_\Omega$  and 
$$J\: G\times T_\Omega\to \Str(V_\C)_0, 
\quad (g,z)\mapsto J(g,z)\:={\partial (g.z)\over \partial z}$$
defines a cocycle (cf.\ [FaKo94, Prop.\ XIII.4.1]). 
Using the simply connectedness of $T_\Omega$, it is not hard to see
that there is a unique lifting of $J$ to a cocycle 
$$ \tilde J\: G\times T_\Omega\to H_\C 
\quad \hbox{ with } \quad \tilde J(\1,z)=\1.$$ 
Let $(\tau, F)$ be a holomorphic representation of $H_\C$. Then we obtain by 
$$J_\tau\: G\times T_\Omega\to \GL(F), \quad 
(g,z)\mapsto J_\tau(g,z)\:=\tau(\tilde J(g,z))$$
a $\GL(F)$-valued cocycle which leads to a representation of 
$G$ on $\Hol(T_\Omega,F)$ given by 
$$(\pi_\tau(g).f)(z)=J_\tau(g^{-1},z)^{-1}.
f(g^{-1}.z) \quad \hbox{ for } \quad g\in G, f\in \Hol(T_\Omega, F).$$ 
For $v \in V$ and $h \in H$ we obtain in particular the simple formulas 
$$(\pi_\tau(v).f)(z)=f(z-v) \quad \hbox{ and } \quad 
(\pi_\tau(h).f)(z)= \tau(h).f(h^{-1}.z). $$

The derived action $d\pi_\tau$ of the Lie algebra $\g$ on $\Hol(T_\Omega,F)$ 
is then given by 
$$ \big(d\pi_\tau(X).f\big)(z) 
= d_1 J_\tau(\1, z)(X).f(z) + \big(\dot\sigma(X).f\big)(z),   $$
where $\dot\sigma(X)(z) := {d\over dt}|_{t=0} \exp(-tX).z$ is the holomorphic 
vector field on $T_\Omega$ corresponding to the action of the one-parameter
subgroup $\exp(\R X) \subeq G$ (cf.\ [Ne99, Prop.\ IV.1.9(ii)]).  

\Lemma III.4. The action of $\g$ on $\Hol(T_\Omega,F)$ is given explicitely by 
polynomial differential operators as follows: 
\item{(i)} $(v.f)(z) = -  df(z)(v)$ for $v \in V \cong \n^+$, 
\item{(ii)} $(X.f)(z) =  d\tau(X).f(z) - df(z)(X.z)$ for $X \in \h$, and 
\item{(iii)} $\big(\theta(v).f\big)(z) =  2 d\tau(z \square v).f(z) 
- df(z)\big(P(z).v\big)$ for $v \in V$, where 
$z \square v \in \str(V_\C) \cong \h_\C$ is defined by 
$a \square b\:= L(ab)+[L(a),L(b)].$

\Proof. In this proof we assume that $v \in V$, $X \in \h$ and $h \in H$. 
According to [FaKo94, p.209], the vector fields corresponding to the
action of $\g$ on $T_\Omega$ are given by 
$$ \dot\sigma(v)(z) = -v, \quad 
\dot\sigma(X)(z) = -X.z \quad \hbox{ and } \quad 
\dot\sigma(\theta.v)(z) = -P(z).v, $$
where the last assertion follows from the fact that the 
inverion $\eta \: T_\Omega \to T_\Omega, \eta(z) = - z^{-1}$ satisfies 
$\eta(g.z) = \theta(g).\eta(z)$ for $g \in G$ and $z \in T_\Omega$. 

We further have 
$$ J_\tau(v,z) = \1, \quad 
J_\tau(h,z) = \tau(h) \quad \hbox{ and } \quad 
J_\tau(\theta.v,z) = \tau\big(\tilde P(v - z^{-1})^{-1} 
\tilde P(z)^{-1}\big). $$
Here the first two formulas are clear from the definitions, and the last 
one follows from $\theta(v).z = \eta(v + \eta(z))$, 
$d\eta(z) = P(z)^{-1}$ ([FaKo94, Prop.\ II.3.3]) and the chain rule. 
This leads to $d_1 J_\tau(\1,z).v = 0$ and $d_1 J_\tau(\1,z).X = d\tau(X)$, 
so that it only remains to verify 
$$ d_1 J_\tau(\1,z).\theta v = 2 d\tau(z \square v). $$
We have 
$$ \eqalign{ d_1 J(\1, z).\theta v  
&= \derat0 P(z^{-1} -  t v)^{-1} P(z)^{-1}  
= -P(z^{-1})^{-1} \derat0 P(z^{-1} - t v) P(z^{-1})^{-1} P(z)^{-1} \cr
&= -P(z) \derat0 P(z^{-1} - t v)  = 2 P(z) P(z^{-1}, v), \cr} $$
so that the assertion follows from Lemma III.5 below. 
\qed

\Lemma III.5. For each $z \in V_\C^\times$ and $u \in V$ we have 
$ P(z) P(z^{-1}, u) = z \square u. $ 

\Proof. In view of [FaKo94, p.\ 147], we have for all $x,y \in V_\C$ and $g \in \Str(V_\C)$
the identities 
$$ P(g.x, g.y) = g P(x,y) g^\top \quad \hbox{ and } \quad 
(g.x \square \theta(g).y) = g(x \square y)g^{-1}. $$
For $z = g.e$ we thus get 
$$ \eqalign{ P(z) P(z^{-1},u) 
&= P(g.e) P\big( (g.e)^{-1}, u\big) 
= g P(e) g^\top  P( \theta(g).e, u) \cr
&= g g^\top \theta(g) P(e, g^\top.u) g^{-1}
= g P(e, g^\top.u) g^{-1}
= g L(g^\top.u) g^{-1}. \cr} $$
On the other hand $e \square x = L(x)$ implies 
$$ z \square u 
= (g.e) \square u 
= g(e \square g^\top.u)g^{-1} 
= g L(g^\top.u) g^{-1} 
= P(z) P(z^{-1},u). $$
Now the assertion follows from $V_\C^\times = \Str(V_\C).e$ 
([FK94, Prop.\ VIII.3.5(i)]). 
\qed 

The main point in the formulas of Lemma III.4 is that they show that 
$d\pi_\tau(\g) \subeq {\cal PD}(V,F)$, acting on $\Hol(T_\Omega,V)$. 
We therefore obtain a representation of $\g$ on $\SO(F)$ by 
$$ X.D := {\cal F}^{-1}(d\pi_\tau(X)).D $$
and further a representation $\rho_\tau$ on $C^\infty(V,F)$ which on 
${\cal S}(V,F)$ is given by 
$$ \la D, \rho_\tau(X).\phi \ra = \la -X.D,\phi \ra 
\quad \hbox{ for } \quad 
\phi \in {\cal S}(V,F), X \in \g. \leqno(3.3) $$

\Remark III.6. For each linear function $\alpha \in V^*$ the subspace 
$e^{\alpha} \Pol(V,F) \subeq C^\infty(V,F)$ is invariant under the 
action of $\g$ because this space is invariant under the whole algebra 
$\PD$. We will see below that the most important space of this type 
arises for $\alpha = - \tr$. 
\qed

The following lemma makes the action of the parabolic subalgebra 
$\n^+ \rtimes \h \subeq \g$ quite explicit. In priciple it is also 
possible to obtain formulas for the action of $\n^- = \theta(\n^+)$, but 
it seems that this can only be done in coordinates. Fortunately, for 
our arguments below it will suffice to have direct access to the action 
of $\h$ and $\n^+$. 

\Lemma III.7. For $\phi \in C^\infty(V,F)$ we have: 
\item{(i)} $\rho_\tau(v).\phi = - i m_v \phi$ for $v \in V \cong \n^+$.  
In particular, $\rho_\tau({\cal U}(\n^+)).\phi= \Pol(V) \cdot\phi.$
\item{(ii)} $\big(\rho_\tau(X).\phi\big)(x) = 
d\tau(\theta(X)).\phi(x) - d\phi(x).(\theta(X).x)$. 

\Proof. (i) We have $d\pi_\tau(v) = - \partial_v$ (Lemma III.4(i)), hence 
${\cal F}^{-1}(d\pi_\tau(v)) = -i m_v$, and therefore 
$$ \la D, \rho_\tau(v).\phi \ra 
=  - \la (-i m_v).D, \phi \ra 
=  \la i m_v.D, \phi \ra 
=  \la D, -i m_v \phi \ra $$
because the pairing $\la \cdot, \cdot \ra$ is sesquilinear. 

\par\nin (ii) Let $v_1, \ldots, v_n$ be a basis of $V$, 
$m_j := m_{v_j}$, and $\partial_j := \partial_{v_j}$. 
Then the differential operator on $\Hol(T_\Omega,F)$ corresponding to $X \in \h$ 
is given by 
$$ (d\pi_\tau(X).f)(z) = 
d\tau(X).f(z) - df(z)(X.z) 
= d\tau(X).f(z) -\Big(\sum_{j,k=1}^n x_{jk} m_k \partial_j.f\Big)(z), $$
where $X.v_j = \sum_{k =1}^n x_{kj} v_k$. The map ${\cal F}^{-1}$ maps this 
operator to 
$$ d\tau(X) -\sum_{j,k=1}^n x_{jk} i \partial_k i m_j 
= d\tau(X) + \sum_{j,k=1}^n x_{jk} \partial_k m_j, $$
and the corresponding action on smooth functions is given by 
$$ \eqalign{ \la D, \rho_\tau(X).\phi \ra 
&= \la \big(-d\tau(X) - \sum_{j,k=1}^n x_{jk} \partial_k m_j\big).D, \phi \ra 
= \la  D, -\tau(X)^*.\phi + \sum_{j,k=1}^n \oline{x_{jk}} m_j \partial_k.\phi \ra \cr
&= \la  D, -\tau(X)^*.\phi + d\phi(\cdot)(X^* \cdot) \ra. \cr} $$
We therefore obtain for $X \in \h$ the formula 
$$ \big(\rho_\tau(X).\phi\big)(x) 
= - d\tau(X)^*.\phi(x) + d\phi(x)(X^*.x) 
= d\tau(\theta(X)).\phi(x) - d\phi(x)(\theta(X).x). 
\qeddis 

From now on we specialize to the case where 
$F := F(\lambda)$ and the representation 
$\tau := \tau_\lambda \: \h \to \End(F(\lambda))$ is 
$\tau_\lambda := \pi_\lambda^\k \circ c_\g$, where 
$(\pi_\lambda^\k, F(\lambda))$ is a unitary highest weight representation 
of $\k_\C$ and $c_\g \in \Aut(\g_\C)$ is the Cayley transform 
(cf.\ Definition III.3(b)). 

In the following we equip $e^{-\tr} \Pol\big(V, F(\lambda)\big)$
with the $\g$-module structure given by the representation 
$\rho_{\tau_\lambda}$ of $\g$ on $C^\infty(V,F(\lambda))$described in (3.3). 

If $\g$ is a Lie algebra, then we define a filtration of ${\cal U}(\g)$ by 
$${\cal U}(\g)^{(n)}=\span_\C\{ X_1\cdot\ldots \cdot 
X_k\: k\leq n, X_i\in \g\}, \quad n \in \N_0.$$
For all $n\in \N_0$ we write $\Pol(V)^{(n)}$ for the space 
of polynomials on $V$ of degree at most $n$. 

\Lemma III.8. There exists a unique isomorphism of $\h$-modules 
$$\Phi\: N(\lambda) \to e^{-\tr} \Pol(V)\otimes F(\lambda) $$
with $\Phi(v)  = e^{-\tr}v$ for $v \in F(\lambda)$. 
Morover, $\Phi$ satisfies 
$\Phi({\cal U}^{(n)}(\h).F(\lambda)) 
=e^{-\tr} \Pol(V)^{(n)}\otimes F(\lambda)$
for all $n\in \N_0$. 

\Proof. To prove the existence of $\Phi$ as a homomorphism of $\h$-modules, in view of 
the description of $N(\lambda)$ as a module induced from the 
$(\h \cap \k)$-module $F(\lambda)$ (Lemma I.6), it suffices to observe that the map 
$$ F(\lambda) \to e^{-\tr} F(\lambda) \subeq C^\infty(V,F), \quad v\mapsto e^{-\tr}.v $$
is a homomorphism of $\h\cap\k$-modules. This is an immediate consequence of the fact that 
the function $e^{-\tr}$ is invariant under $H\cap K \subeq \Aut(V)$ because 
$(\h \cap \k).e = \{0\}$. 

Next we show that $\Phi$ is injective. 
Let $v_1, \ldots, v_m$ denote a basis of $V$, so that the 
Jordan left multiplications $L(v_j)$ form a basis of $\h\cap \p$. And write 
${\cal U}(\h \cap \p)$ for the subspace of ${\cal U}(\h)$ spanned by the ordered 
products of the form $L(v_1)^{k_1} \cdots L(v_m)^{k_m}$. We will use the vector space 
isomorphism 
$$ N(\lambda)
\cong {\cal U}(\h)\otimes_{{\cal U}(\h\cap\k)} F(\lambda) 
\cong {\cal U}(\h\cap \p)\otimes F(\lambda). $$
For $X\in \h$ and $p\in \Pol(V)\otimes F(\lambda)$ we obtain from Lemma II.7(i)
for all $v \in V$ the formula 
$$\eqalign{
&\ \ \ \ d\pi_\lambda(X).(e^{-\tr} p)(y)\cr
&=d\tau_\lambda(\theta(X)). e^{-\tr(y)} p(y)
-(y,X.e) e^{-\tr(y)} p(y) -e^{-\tr(y)} dp(y).(\theta(X).y).\cr}\leqno(3.4)$$
Inductively this leads for $u=L(x_1)\cdot \ldots \cdot L(x_n)$, $x_j \in V$, to 
$$\Phi(u\otimes v)(y)= e^{-\tr(y)}\Big(\Big(\prod_{j=1}^n (y,x_j)\Big) v + \hbox{polynomials 
of order $<n$}\Big).\leqno(3.5)$$
From (3.5) we deduce that $\Phi$ is injective. 

We also conclude inductively with (3.5) that 
$$\Phi({\cal U}^{(n)}(\h)\otimes F(\lambda)) =e^{-\tr} \Pol(V)^{(n)}\otimes F(\lambda)$$
for all $n\in \N_0$, and this implies in particular that $\Phi$ is surjective.
\qed 

\Remark III.9. There are several actions of the group $H$ which will be used in the following. 

\par\nin (a) First we recall the action on $\Hol(T_\Omega, F(\lambda))$ given by 
$$ (\pi_\tau(h).f)(z) = \tau_\lambda(h).f(h^{-1}.z), $$
and on the space $C^\infty(V,F(\lambda))$ we consider the action given by 
$$ (h.\phi)(x) = \tau_\lambda(\theta(h)).f(\theta(h)^{-1}.x) $$
which is the integrated form of the representation $\rho_\tau$ of $\h$ (Lemma III.7). 
On the space $C^\infty(V)$ we consider the action given by 
$(h.\phi)(x) := \phi(h^{-1}.x)$. 
Via $(h.D)(\phi) := D(h^{-1}.\phi)$ we further obtain an action on ${\cal D}'(V,F(\lambda))$ which 
is compatible with the action of $\h$ and given on 
$D \otimes v$ by 
$$ \eqalign{ \la h.(D \otimes v), \phi \otimes w \ra 
&= \la D \otimes v, h^{-1}.(\phi \otimes w) \ra 
= \la D \otimes v, \theta(h)^{-1}.\phi \otimes \tau_\lambda(h)^*.w \ra \cr
&= \la \theta(h).D \otimes \tau_\lambda(h).v, \phi \otimes w \ra, \cr} $$
i.e., 
$$ h.(D \otimes w) = \theta(h).D \otimes \tau_\lambda(h).v. $$
This action preserves the subspace $\SO(F(\lambda))$, and the Fourier transform 
$${\cal F} \: \SO(F(\lambda)) \to \Hol(T_\Omega, F(\lambda))$$ 
intertwines the $H$-actions on
both spaces. On the level of the derived representation this is satisfied by definition, and 
on the group level it either follows from the connectedness of $H$ and the smoothness of the 
actions, or directly by calculation. 

\par\nin (b) Further we consider on the space 
$\End(F(\lambda))$ the $H$-action given by 
$h.A := \tau_\lambda(h)A\tau_\lambda(h)^*,$ and the induced action on 
$\SO(\End(F(\lambda))$ given by 
$$ h.(D \otimes A) = \theta(h).D \otimes h.A. $$

\par\nin (c) We also have a natural map 
$$ \SO\big(\End(F(\lambda))\big) 
\times C^\infty(V,F(\lambda)) \to {\cal D}'(V,F(\lambda)), \quad 
(D \otimes A).(\phi \otimes v) \mapsto \phi D \otimes A.v. $$
The observation that 
$h.(\phi D) = (h.\phi) (h.D)$ holds for the multiplication of distributions with 
smooth functions,  
together with (a) and (b), shows that this map is $H$-equivariant. 
\qed

\Theorem III.10. There exists a unique 
tempered distribution $R_\lambda \in \SO\big(\End(F(\lambda)\big)$ with 
$$ {\cal F}(R_\lambda)(z) = \tau_\lambda\big(\tilde P(-i z)\big) \quad \hbox{ for all } \quad 
z \in T_\Omega. $$
This distribution is $H$-invariant. 

\Proof. The existence of $R_\lambda$ follows from [HiNe99, Th.\ V.8]. 
The formula $P(h.x) = h P(x) h^\top$ for $h \in \Str(V)$ and $x \in V$ implies that 
$$ {\cal F}(R_\lambda)(z)
= \tau_\lambda(h) {\cal F}(R_\lambda)(h^{-1}.z) \tau_\lambda(h)^* \quad \hbox{ for } \quad 
z \in T_\Omega, h \in H. $$
It is easy to see that the action of $H$ on $\Hol(T_\Omega,\End(F(\lambda)))$ given by the 
formula on the right hand side is intertwined by the Fourier transform with the 
action on $\SO\big(\End(F(\lambda))\big)$, so that the injectivity of the Fourier 
transform implies that $R_\lambda$ is fixed by $H$. 
\qed

\Theorem III.11. {\rm(Equivariance Theorem)} The map 
$$ \Psi \: N(\lambda) \to \SO(F(\lambda)), \quad 
v \mapsto R_\lambda.\Phi(v) $$ 
is $\g$-equivariant, and its range 
$R_\lambda.e^{-\tr}\Pol(V,F(\lambda))$ is isomorphic to $L(\lambda)$. 

\Proof. Since the Fourier transform ${\cal F} \: \SO(F(\lambda)) \to \Hol(T_\Omega,F(\lambda))$
is by definition equivariant with respect to the action of $\g$ on both sides, it suffices 
to show that the composition 
$$ {\cal F} \circ \Psi \:  N(\lambda) \to \Hol(T_\Omega, F(\lambda)) $$
is equivariant. The operator 
$$ C_\lambda\: \Hol({\cal D}, F(\lambda))\to \Hol(T_\Omega, F(\lambda)), \quad 
(C_\lambda.f)(z)\:= \tau_\lambda\Big(\tilde P\big({z+ie\over 2i}
\big)\Big)).f ((z-ie)(z+ie)^{-1})$$
intertwines the natural representations of 
the covering group of $\Aut({\cal D})_0$, resp.\ $G$,  on both spaces 
(cf.\ [FaKo94, Lemma X.4.4]). 
The subspace $F(\lambda)$, as a $\tilde \k$-module 
is embedded into $\Hol({\cal D}, F(\lambda))$ as the subspace 
of the constant functions, the lowest $\tilde \k$-type, and this subspace 
generates a simple highest weight module of $\tilde\g$ of highest weight 
$\lambda$ (cf.\ [Ne99, Prop.\ XII.2.2]). 
Let $d_\lambda > 0$ be the number determined by 
$\tau_\lambda({1\over 4}\1) = d_\lambda \1$. Then we obtain for 
$v \in F(\lambda) \subeq \Hol({\cal D},F(\lambda))$ the formula 
$$ \eqalign{ C_\lambda(v)(z)
&=\tau_\lambda\Big(\tilde P\Big({z+ie\over 2i}\Big)\Big).v
=\tau_\lambda({\textstyle{1\over 4}}\1) \tau_\lambda\Big(\tilde P\Big({z+ie\over i}\Big)\Big).v\cr
&=d_\lambda {\cal F}(R_\lambda)(z+ie).v 
=d_\lambda {\cal F}(R_\lambda.v)(z+ie)
=d_\lambda {\cal F}(R_\lambda.e^{-\tr} v)(z).\cr} $$
This implies that the lowest $K$-type for the $G$-representation in 
$\Hol(T_\Omega,F(\lambda))$ is given by the subspace 
${\cal F} \circ \Psi\big(F(\lambda))$
which is annihilated by $\p^+$. 

Now the definition of $N(\lambda)$ as a $\g$-module induced from the $\p^+ \rtimes \k_\C$-module 
$F(\lambda)$ implies the existence of a $\g$-equivariant map 
$\Psi_1 \: N(\lambda) \to \Hol(T_\Omega, F(\lambda))$. Its range is the $\g$-module 
generated by ${\cal F} \circ \Psi\big(F(\lambda))$ which is a simple 
highest weight module of highest weight $\lambda$, hence isomorphic to $L(\lambda)$. 

Since $R_\lambda$ is $H$-invariant (Theorem III.10), the multiplication map 
$$ e^{-\tr}\Pol(V,F(\lambda)) \to \SO(F(\lambda)), \quad 
\phi \mapsto R_\lambda.\phi $$
is $H$-equivariant (cf.\ Remark III.9(c)). Now Lemma III.8 implies that 
$\Psi$ is $\h$-equivariant. In view of Proposition I.5, this shows that 
${\cal F} \circ \Psi$ maps $N(\lambda) = {\cal U}(\h).F(\lambda)$ surjectively onto 
${\cal U}(\h).({\cal F} \circ \Psi)(F(\lambda)) \cong L(\lambda)$. 
Therefore ${\cal F} \circ \Psi$ and $\Psi_1$ can be viewed as two $\h$-equivariant linear maps 
$N(\lambda) \to F(\lambda)$ which coincide on $F(\lambda)$. 
Now $\ker({\cal F} \circ \Psi - \Psi_1)$ is an $\h$-invariant subspace of $N(\lambda)$ containing 
$F(\lambda)$, hence coincides with $N(\lambda)$ (Proposition I.5). This proves 
that ${\cal F} \circ \Psi = \Phi_1$ is $\g$-equivariant. 
\qed

The Equivariance Theorem has several interesting consequences. 

\Corollary III.12. {\rm(a)} The subspace 
$R_\lambda. e^{-\tr}\Pol(V,F(\lambda))$ is a $\g$-submodule of 
$\SO(F(\lambda))$ isomorphic to $L(\lambda)$. 

\par\nin {\rm(b)} $\ker \Psi = J(\lambda)$ is the maximal submodule of $N(\lambda)$. 

\par\nin {\rm(c)} $J(\lambda) = \{ \phi \in N(\lambda) \cong e^{-\tr} \Pol(V,F(\lambda)) 
\: R_\lambda.\phi = 0\}$. 
\qed

\Lemma III.13. Let $f=\sum_{j=0}^nf_j\in J(\lambda)$ with 
$f_j=e^{-\tr}p_j$ and $p_j\in \Pol(V, F(\lambda))$ homogeneous 
of degree $j$. Then we have $f_j\in J(\lambda)$ for all $0\leq j\leq n$. 

\Proof. Recall that $L(e)$ is the generator of $\z(\h)\subeq \h\cap\p$. 
From (3.4) we obtain for all $x\in \supp R_{\lambda}$: 
$$ \eqalign{ d\rho_\lambda(L(e).f)(x)
&= -d\tau_\lambda(L(e)).f(x) -(x,e) f(x) +\sum_{j=0}^n j f_j(x).\cr} $$
Clearly we have $d\tau_\lambda(L(e)).f(x) \in \C f(x) \subeq 
J(\lambda)$, and by Corollary III.12(c) we also 
have $\tr \cdot f\in J(\lambda)$, so that 
$\sum_{j=0}^n j f_j\in J(\lambda)$. This means that 
$J(\lambda)$ is invariant under the operator 
$E(f) := \sum_{j=0}^n j f_j$, hence adapted to the eigenspace decomposition of $E$. 
\qed

\sectionheadline{IV. Unitary highest weight representations for Cayley type spaces}

\nin In this section we now turn to unitary highest weight representations 
in the setup studied in Section II. We will recall the realization of 
unitary highest weight representations in spaces of vector valued 
square integrable functions on a homogeneous cone.  
Then we give a new and quite explicit characteriztion of the $K$-finite vectors 
in the cone realization which will easily follow from the more general 
results in Section II. 

We keep the notation from Section III. In particular 
$G$ is the simply connected covering group of $\Aut(T_\Omega)_0$ and 
$K$ and $H$ are connected subgroups of $G$. 

Let $(\pi_\lambda, {\cal H}_\lambda)$ be a unitary highest weight 
representation of $G$ with highest weight $\lambda\in i\t^*$ and with 
respect to $\hat\Delta_k^+$. We write $(\pi_\lambda^K, F(\lambda))$ for the corresponding 
irreducible $K_\C$-representation with highest weight $\lambda$. 
Thus Definition III.3(b) 
implies that the prescription $d\tau_\lambda \:=d\pi_\lambda^K\circ c_\g^{-1}$ integrates 
to a holomorphic representation $(\tau_\lambda, F(\lambda))$ of $H_\C$.

\Theorem IV.1. Let $R_\lambda \in \SO(\End(F(\lambda))$ denote the distribution with 
$$ {\cal L}(R_\lambda)(z) = \tau_\lambda(\tilde P(-iz)), \quad z \in T_\Omega. $$
Then the highest weight module $L(\lambda)$ of $\g$ is unitarizable if and only if 
$R_\lambda$ is a measure with values in the cone $\Herm^+(F(\lambda))$ of positive operators on 
$F(\lambda)$. If this is the case, then the following assertions hold: 
\item{(i)} Let  
$$L^2(\oline \Omega, R_\lambda)\:=
\Big\{ f\: \oline \Omega\to F(\lambda) \hbox{\rm\ meas.}\: 
\int_{\oline\Omega} \la dR_\lambda(x).f(x), f(x)\ra<\infty\Big\}. $$
Then the  Fourier transform 
$${\cal F}_\lambda\: L^2(\oline \Omega, R_{\lambda})\to {\cal
H}_\lambda\subeq \Hol(T_\Omega, F(\lambda)), \quad 
{\cal F}_\lambda(f)(w)\:=\int_{\oline\Omega} e^{i(w,x)} dR_{\lambda}(x).f(x) $$
yields an isomorphism on the reproducing kernel Hilbert space ${\cal H}_\lambda$ with kernel 
$$ \tau_\lambda\Big(\Big(\tilde P({z- \oline w\over 2i}\Big)\Big). $$
\item{(ii)} If we realize the unitary highest weight representation $(\pi_\lambda, {\cal H}_\lambda)$ of 
$G$ in the space $L^2(\oline\Omega, R_\lambda)$ as in {\rm (i)}, then the 
action of the parabolic subgroup $V \rtimes H$ is given by 
$$ (\pi_\lambda(v,h).f)(x) = e^{-i(v,x)} \tau_\lambda(\theta(h)).f(\theta(h)^{-1}.x). $$

\Proof. For the first statement we have to combine several results. 
First [Ne99, Th.\ XII.2.6] says that $L(\lambda)$ is unitarizable if the 
natural $\End(F(\lambda))$-valued kernel on ${\cal D}$ is positive definite, which 
by the Cayley transform is equivalent to the positive definiteness of the function 
$\tau_\lambda \circ \tilde P$ on $T_\Omega$. This in turn is equivalent to 
$R_\lambda$ being a $\Herm^+(F(\lambda))$-valued measure 
([HiNe99, Th.\ V.12], see also [Cl95, Th.\ 3.3]). 

\par\nin (i) [Cl95, p.\ 233] 
\par\nin (ii) This is immediate from [Cl95, Th.\ 3.4].\qed 

The preceding result is the appropriate version of Clerc's results that we need in our 
context. For the scalar case ($\dim F(\lambda) = 1$), similar results were obtained by 
Vergne and Rossi in [VR76]. 

\Remark IV.2. Suppose that $R_\lambda$ is a $\Herm^+(F(\lambda))$-valued measure. 
Then for each element $\phi \in L^2(\oline\Omega,R_\lambda)$ we can view 
$R_\lambda.\phi$ as an element of $\SO(F(\lambda))$ and obtain 
$$ {\cal F}(R_\lambda.\phi) = {\cal F}_\lambda(\phi). $$
Since the maps ${\cal F}$ and ${\cal F}_\lambda$ are $\g$-equivariant, it follows 
from Lemma III.7 and Theorem III.11 that the action of $\g$ on the subspace 
$L^2(\oline\Omega, R_\lambda)^\infty$ of smooth vectors is given by the 
natural action on $\SO(F(\lambda))$. For the subalgebra $\n^+ \rtimes \h$ it is given by 
the formulas in Lemma III.7. 
\qed

In the following lemma the assumption that $L(\lambda)$ is unitarizable is used to 
ensure that $J(\lambda)$ is a cyclic module. 

\Lemma IV.3. Suppose that $N(\lambda)$ is unitarizable and of reduction level 
$k$. Then 
there exists an element $v_k\in {\cal U}(\h)^{(k)}.v_\lambda\bs 
{\cal U}(\h)^{(k-1)}.v_\lambda $ such that the maximal 
submodule $J(\lambda) \subeq N(\lambda)$ is given by 
$J(\lambda)={\cal U}(\h).v_k.$

\Proof. In view of Lemma III.8, we may identify in the following the $\h$-module 
$N(\lambda)$ with the space $e^{-\tr} \Pol(V)\otimes F(\lambda)$. 
Comparing dimensions in the following chain of inclusions 
$$ {\cal U}(\p^-)^{(n)} \otimes F(\lambda) 
= {\cal U}(\g)^{(n)}.F(\lambda) 
\supeq {\cal U}(\h)^{(n)}.F(\lambda) 
= e^{-\tr} \Pol(V)^{(n)} \otimes F(\lambda) $$
yields equality, so that the assertion follows from 
Proposition I.5 and Theorem II.4.
\qed 

\subheadline{Identification on central half lines}

Let $\lambda^0$ be the highest weight of $\pi^\k_{\lambda^0}$ with respect to 
the Cartan subalgebra $\t$. 
Then 
$$\lambda^0 \circ c_\g\res_\b = \sum_{j=1}^r {m_j\over 2} \eps_j \quad \hbox{ with } \quad 
m_j - m_i \in \N_0, \quad i < j \leqno(4.1) $$
(Definition III.3(b)). 
We assume, in addition, that $m_r = 0$. 
Then $m_j \in -\N_0$ for all $j$. For 
$u\in \R$ we then set 
$$\lambda\:=\lambda^u\:=\lambda^0 - u\zeta,$$
where $\zeta \in [\k,\k]^\bot$ is determined by $\zeta\circ c_\g\res_\b = {1\over 2}\sum_j \eps_j$ 
which is equivalent to $\zeta(\check \gamma)= 1$ for $\gamma \in \hat\Delta_n^+$ long because 
for the highest root $\gamma \in \hat\Delta_n^+$ the corresponding restricted root 
is $\eps_r$ with the coroots $2L(c_r)$ (cf.\ [Kr99b, Sect.\ IV]). 
Then $\lambda$ is the highest weight of a 
finite-dimensional irreducible representation $(\tau_\lambda, 
F(\lambda))$ of $H$ which can be written as 
$$\tau_\lambda=\tau_{\lambda^0}\otimes\tau_{-{u\over 2}\eps},$$ where 
$\tau_{-{u\over 2}\eps}$ is a one-dimensional representation. 
Accordingly, we may identify $F(\lambda^u)$ for all $u\in \C$ with 
the fixed vector space $F_0\:=F(\lambda^0)$.

\par  Recall the definition of the Jordan algebra determinant 
$\Delta\:=\det$ from [FaKo94, p.\ 29]. 
For each $u\in \R $ we denote by $R_u$ the Riesz distribution
with parameter $u \in \C$ on $\oline \Omega$ which are uniquely determined by 
$$ {\cal F}(R_u)(z) = \Delta(-iz)^{-u}, \quad z \in T_\Omega, u \in \C $$
(cf.\ [FaKo94, Sect.\ VII.2]). 
By a theorem of Gindikin (cf.\ [FaKo94, Th.\ VII.3.2]), $R_u$ is a
positive measure if and only if 
$$u\in \Big\{ 0, {{d\over 2}}, \ldots (r-1){{d\over 2}}\Big\}\cup](r-1){{d\over 2}}, \infty[, $$ 
where $d\:=\dim\g^{{1\over 2}(\eps_i+\eps_j)}$, $i\neq j$, is obtained from the 
restricted root decomposition of $\g$ with respect to $\b$. Moreover we have 
$\supp(R_u)=\oline\Omega$ with $R_u(\oline\Omega\bs \Omega)=0$ for 
$u>(r-1){d\over 2}$ and $\supp(R_u)=\oline {{\cal O}_k}$ with   $R_u(\oline{{\cal O}_k}\bs 
{\cal O}_k)=0$ for $u=k{d\over 2}$, $0\leq k\leq r-1$. 

\par By a more detailed study of the vector valued measure $R_{\lambda}$, 
one can retrieve useful information on $R_{\lambda}$ 
which we collect in the next theorem. 

\Theorem IV.4. If $R_{\lambda^u}$ is a $\Herm^+(F({\lambda^u}))$-valued measure, then: 
\item{(i)} There exists a polynomial map 
$D_u \: V \to \Herm^+(F({\lambda^u}))$ with $R_\lambda = D_u R_u$. 
\item{(ii)} $u\in \{ 0, {d\over 2},\ldots, (r-1){d\over 2}\} \cup ] (r-1) {d\over 2}, \infty[$.
\item{(iii)} $\supp(R_{\lambda^u})=\supp(R_{u})=\cases
{\oline \Omega  & for $u>(r-1) {d\over 2} $\cr 
\oline{{\cal O}_k} & for $u=k {d\over 2}, k \in \{0,\ldots, r-1\}$ \cr}.$

\Proof. [HiNe99, Th.\ V.17]
\qed 

\subheadline{The space of $K$-finite vectors}

We now determine explicitely the space of $K$-finite vectors for a  
unitary highest weight representations of $G$ realized in $L^2(\oline\Omega, R_\lambda)$. 
Since the 
$K$-action in the cone realization is rather nasty (cf.\ [DG93]), there is no obvious 
description of this space. Our approach is based on the action of the maximal parabolic 
subalgebra $\h\ltimes \n^+$ which led in Section II already to the right subspace of 
$\SO(F(\lambda))$. 

\par Recall the decomposition of the ring $\Pol(V)$ of polynomials into 
irreducible $H$-modules (cf.\ [FaKo90, Th.\ 2.1], [FaKo94, Th.\ XI.2.4]): 
The $H$-decomposition is muliplicity free and given by 
$$\Pol(V)=\bigoplus_{{\bf m}\geq 0}  \Pol(V)_{\bf m} \leqno(4.2) $$
where ${\bf m} = (m_1, \ldots, m_r)$ and the condition ${\bf m} \geq 0$ means that 
$m_1\geq m_2\geq \ldots\geq m_r\geq 0$ and $m_i\in \Z$. 
The space $\Pol(V)_{\bf m}$ is a real 
highest weight module of $H$ with highest weight 
$\lambda_{\bf m}=\sum_{j=1}^r m_j \eps_j \in \b^*$  w.r.t. $-\Sigma^+$ 
which is generated by the highest weight vector 
$\Delta_{\bf m}$ (for the definition of the generalized power functions $\Delta_{\bf m}$ 
we refer to [FaKo94, p.\ 122]).   
\par Recall that the Jordan algebra trace $\tr$ on $V$ satisfies 
$\tr(x)=(x,e)$ for all $x\in V$.  

\Theorem  IV.5. {\rm (Description of $K$-finite vectors)}
Suppose that $(\pi_\lambda, {\cal H}_\lambda)$ is 
a unitary highest weight representation of the simply 
connected hermitian group $G$ of tube type realized as $L^2(\oline\Omega, R_{\lambda})$. 
Then the space $L^2(\oline\Omega, R_\lambda)^K$ of  $K$-finite vectors of $(\pi_\lambda, 
{\cal H}_\lambda)$  is the image of the space $e^{-\tr} \Pol(V, F(\lambda))$
in $L^2(\oline\Omega, R_{\lambda})$, and 
a highest weight vector of $(\pi_\lambda, 
{\cal H}_\lambda)$ is given by the function $e^{-\tr}.v_\lambda,$
where $v_\lambda$ is a highest weight vector of $(\tau_\lambda, F(\lambda))$.

\Proof. We identify $L^2(\oline\Omega, R_\lambda)$ with the subspace 
$$ R_\lambda.L^2(\oline\Omega, R_\lambda) \subeq \SO(F(\lambda)) $$
which is compatible with the $\g$-action on the space of 
smooth vectors, because the Fourier transform 
$$ {\cal F} \: R_\lambda.L^2(\oline\Omega, R_\lambda)^\infty \to 
\Hol(T_\Omega, F(\lambda)) $$
is $\g$-equivariant (cf.\ Remark IV.2). 

Since $R_\lambda$ is a tempered distribution, all restrictions of the functions in 
$e^{-\tr} \Pol(V, F(\lambda))$ are contained in $L^2(\oline\Omega, R_\lambda)$. 
In view of Theorem III.11, the subspace 
$R_\lambda.e^{-\tr}\Pol(V,F(\lambda))$ is a $\g$-module isomorphic to 
$L(\lambda)$, hence contained in $L^2(\oline\Omega,R_\lambda)^\infty$, and therefore in 
$L^2(\oline\Omega,R_\lambda)^K$. 
On the other hand, the fact that the representation of $G$ on 
$L^2(\oline\Omega,R_\lambda)$ is a unitary highest weight representation implies that
$L^2(\oline\Omega,R_\lambda)^K$ is a simple $\g$-module, which yields equality. 
The description of the highest weight vector follows from 
Theorem III.11. 
\qed

In the scalar case Theorem IV.5 has been proved by completely different methods 
by J.~Faraut and A.~Kor\'anyi (cf.\ [FaKo94, Prop.\ XIII.3.2].) 

For every $0\leq k\leq r-1$ we write 
$$ I_k := \bigoplus_{{\bf m}\geq 0, m_{k+1}>0}\Pol(V)_{\bf m} = \{ f \in\Pol(V)\: 
f\res_{{\cal O}_k} = 0\} $$
for the ideal defined by ${\cal O}_k$. 
For the asserted equality, we note that the ideal defined by ${\cal O}_k$ 
is $H$-invariant, hence adapted to the $H$-decomposition (4.2) of $\Pol(V)$. 
Therefore one only has to note that for ${\bf m} \geq 0$ the function 
$\Delta_{{\bf m}}$ vanishes on ${\cal O}_k$ if and only if $m_{k+1} > 0$ 
(cf.\ [HiNe99, Prop.\ I.4]). 

The following proposition provides a description of the maximal submodule $J(\lambda) \subeq 
N(\lambda)$ in the picture of Theorem IV.5. 

\Proposition IV.6. Suppose that $L(\lambda)$ is unitarizable and write 
$R_\lambda = D_u R_u$ with the polynomial density function $D_u$. 
Identifying $N(\lambda)$ with $e^{-\tr}\Pol(V,F(\lambda))$, we have 
$$J(\lambda)=\{ e^{-\tr} f\in N(\lambda)\: D_u.f =0\  \hbox {on $\supp(R_u)$}\}. $$
For $\lambda = \lambda^u$, $u=k{d\over 2}$, $k \in \{0,\ldots, r-1\}$, we obtain in particular 
$J(\lambda) \supeq e^{-\tr} I_k \otimes F(\lambda).$
 
\Proof. In view of Theorem IV.5, the natural $\g$-contravariant hermitian form on 
$$ N(\lambda)\cong e^{-\tr} \Pol(V, F(\lambda)\big)$$ 
is given by 
$$ \la e^{-\tr} f,e^{-\tr}g \ra 
= \int_{\oline\Omega} e^{-2 \tr(x)} \la dR_\lambda(x).f(x), g(x)\ra
= \int_{\oline\Omega} e^{-2 \tr(x)} \la D_u(x).f(x), g(x) \ra\, dR_u(x).  $$
This implies the assertion because 
$J(\lambda) = \{ e^{-\tr} f \in N(\lambda) 
\: \la e^{-\tr} f,e^{-\tr}f \ra = 0\}.$

For the special case described in the second part of the statement, we only have to 
observe that $\supp(R_\lambda) = \oline{{\cal O}_k}$ which implies that 
$R_\lambda.e^{-\tr} I_k\otimes F(\lambda) = \{0\}$
in this case. This completes the proof. 
\qed

\sectionheadline{V. The classification for Cayley type spaces}

\nin In this section we give an explicit description of the $H$-fixed distribution vector of 
a spherical unitary highest weight representation $(\pi_\lambda, {\cal H}_\lambda)$ of a 
simply connected hermitian group $G$ of tube type  
(cf.\ Theorem V.1). In the case where $\supp(R_\lambda)$ is smaller than 
$\oline\Omega$ our approach yields a complete classification 
of all spherical unitary highest weight representations of $G$. 

\subheadline{The $H$-invariant distributions} 

In this section we assume that $(\tau_\lambda, F(\lambda))$ is a $H\cap K$-spherical 
representation of $H$. In view of the Appendix A in [HiNe99], this means 
that it can be viewed as a real highest weight representation 
with respected to the restricted root decomposition of $\h$ with respect to $\b$. 
We write 
$$\lambda = \lambda^s := \sum_{j=1}^r {m_j-s\over 2} \eps_j \leqno(5.1) $$ 
with $m_1 \leq \ldots \leq m_k < m_{k+1} = \ldots = m_r = 0$ and $s \in \R$. 
Then the corresponding representation satisfies 
$\tau_{\lambda^s} = \tau_{\lambda^0} \otimes \det^{-s{r \over 2\dim V}}$
(cf.\ [HiNe99, Sect.\ I]), so that we may identify all the spaces 
$F(\lambda^s)$ with $F_0 := F(\lambda^0)$. 

According to [KN\'O97, Th.\ II.11], there exists for each 
$v_0 \in F(\lambda^s)^{\h \cap \k} \cong F_0^{\h \cap \k}$
a unique $H$-invariant holomorphic function 
$$ f_s \: T_\Omega \to F(\lambda^s) \quad \hbox{ with } \quad 
f_s(ie) =v_0. $$
Note that the $H$-invariance of $f$ with respect 
to the action $(h.f)(z) = \tau_{\lambda^s}(h).f(h^{-1}.z)$ means that 
$f$ is an $H$-equivariant map. The uniqueness of $f$ is a direct consequence of the 
fact that $f$ is uniquely determined by its restriction to $i\Omega = H.ie \subeq T_\Omega$. 

\Theorem V.1. There exists a weakly holomorphic map 
$\C \to {\cal S}'\big(V,F(\lambda^0)\big),s \mapsto R_{\lambda^s}^H$ 
with ${\cal F}(R_{\lambda^s}^H)\res_{T_\Omega} = f_s$ 
and $R_{\lambda^s}^H$ is supported by $\oline\Omega$. 
Suppose that 
$$ m_1 \leq \ldots \leq m_k < m_{k+1} = \ldots = m_r =0. $$
Then the tempered distribution $R_{\lambda^s}^H$ is a measure if 
and only if 
$$ {s\over 2} \in \{k{\textstyle{d\over 2}}, \ldots, (r-1){\textstyle{d\over 2}}\}
\cup \big](r-1){\textstyle{d\over 2}}, \infty\big[. $$
If this is the case, then 
$$ R_{\lambda^s}^H = ({\textstyle{s\over 2}})_{-{\bf m}}^{-1} \eta \cdot R_{s\over 2}, $$
where 
$$ ({\textstyle{s\over 2}})_{-{\bf m}} 
= \prod_{j = 1}^r \prod_{l=0}^{-m_j - 1} \Big({s\over 2}-(j-1){d\over 2}+l\Big)
= \prod_{j = 1}^k \prod_{l=0}^{-m_j - 1} \Big({s\over 2}-(j-1){d\over 2}+l\Big), $$
$\eta \: V \to F(\lambda^0)$ is a polynomial map with 
$\eta(g.e) = \tau_{\lambda^0}(\theta(g)).v_0$ for $g \in H$, and 
$R_s$, $s \in \C$, are the Riesz measures on $\oline\Omega$ with 
${\cal F}(R_s)(z) = \Delta(-i z)^{-s}.$

\Proof. The existence of the weakly holomorphic family $R_{\lambda^s}^H$, 
$s \in \C$, follows from [HiNe99, Prop.\ V.2]. The characterization of those 
parameters for which $R_{\lambda^s}^H$ is a measure is a consequence of 
[HiNe99, Props.\ III.8 and V.5]. That the density 
$\eta$ is a polynomial follows from [HiNe99, Props.\ II.8]. 
\qed

\Remark V.2. (a) The polynomial $({s\over 2})_{-{\bf m}}$ vanishes in the points of the 
form $(j-1){d\over 2} - l$, where 
$1 \leq j \leq k$ and $0 \leq l \leq - m_j - 1$. Therefore the maximal 
zero is $s = (k-1)d$, which is the place, where the 
distribution $R_{\lambda^s}^H$ no longer is a measure. 
The values 
$s = jd$ correspond to measures on the boundary components 
${\cal O}_j$, $k \leq j \leq r-1$. For $j < k$ the density function 
$\eta$ vanishes on ${\cal O}_j$ ([HiNe99, Prop.\ II.12]), which 
somehow cancels with the zero of $({s\over 2})_{-{\bf m}}$, but the resulting 
distribution is not a measure.  

\par\nin (b) Since the Fourier transform 
${\cal F} \: {\cal S}'_{\oline\Omega}(V,F(\lambda)) \to \Hol(T_\Omega, F(\lambda))$
is injective and equivariant with respect to the $H$-actions on both sides (Remark III.9), 
the uniqueness of an $H$-invariant holomorphic function implies that 
$$  \dim {\cal S}'_{\oline\Omega}(V,F(\lambda))^H = 1.
\qeddis 

Let $F(\lambda)$ be a spherical $H$-representation with highest weight $\lambda$ and 
$0 \not= v_0 \in F(\lambda)^{\h \cap \k}$. 
We have seen in the preceding section that there exists a unique 
$H$-invariant tempered distribution $R_\lambda^H \in {\cal S}'_{\oline\Omega}(V, F(\lambda))$ with 
${\cal F}(R_\lambda^H)(ie) = v_0$.

\Lemma V.3. The prescription 
$$  \nu \: N(\lambda) \cong e^{-\tr} \Pol(V, F(\lambda)) \to \C, \quad 
e^{-\tr} f \mapsto 
\la R_\lambda^H, e^{-\tr} f  \ra $$is a non-zero element of 
$\big(N(\lambda)^\sharp\big)^\h$. 

\Proof. It follows from [HiNe99, Lemma B.1] that $\nu$ is well 
defined because $e^{-\tr}\res_{\oline\Omega}$ extends to a Schwartz function on $V$. 
Finally $\nu$ is annihilitad by $\h$ since $R_\lambda^H$ is $H$-fixed. \qed 

\Corollary V.4. $L(\lambda)$ is spherical if and only if 
$\la R_\lambda^H, J(\lambda) \ra = \{0\}.$
\qed

\subheadline{The classification of unitary highest weight modules with spherical 
$F(\lambda)$} 

Let $\gamma_1, \ldots, \gamma_r \in \hat\Delta_n^+$ be a maximal system of 
orthogonal roots, where $\gamma_r$ is the highest root and 
$\gamma_{j-1}$ is maximal among all positive non-compact roots orthogonal 
to $\gamma_j, \ldots, \gamma_r$. Then a linear functional 
$\lambda \in i \t^*$ is invariant under $(-\tau)$ if and only if 
$\lambda \in \span \{\gamma_j \: j =1,\ldots, r\}$ because the real rank of 
$\g$ equals $r$ (cf.\ Definition III.3(b)). Therefore it is of the form 
$$ \lambda = \sum_{j=1}^r {\lambda(\check \gamma_j)\over 2} \gamma_j. $$
Since $\gamma_j \circ c_\g\res_\b = \eps_j \in \b^*$, we have 
$$ \lambda_\b := \lambda \circ c_\g\res_\b = 
\sum_{j=1}^r {\lambda(\check \gamma_j)\over 2} \eps_j. $$

That $F(\lambda)$ is a $\h \cap \k$-spherical highest weight module of 
$\k$ is equivalent to $-\tau.\lambda = \lambda$ and the condition that 
$\lambda_\b \in \b^*$ is the highest weight of a spherical representation 
of $\h$ which is equivalent to 
$$ \lambda(\check \gamma_j) - \lambda(\check \gamma_i) \in 2 \N_0  
\quad \hbox{ for } \quad i < j \leqno(5.2) $$
(cf.\ Definition III.3(a)). 
For $\lambda = \lambda^u$ as in (5.1) this means that 
$m_j = \lambda(\check \gamma_j)$. 

The following theorem shows that the classification of unitary highest weight modules 
with spherical $F(\lambda)$ is much simpler than the general case.  

\Theorem V.5. If $F(\lambda^0)$ is $\h \cap \k$-spherical, 
$\lambda_\b = \sum_{j=1}^r {m_j \over 2} \eps_j$, and 
$k \in \{0,\ldots, r-1\}$ is determined by 
$m_1 \leq \ldots \leq m_k < m_{k+1} = \ldots = m_r,$
then 
$$l(\lambda^0)
=\Big\{ k d, (2k+1){d\over 2}, \ldots, (k + r-1){d\over 2}\Big\} 
\dot\cup \Big](k + r-1){d\over 2}, \infty \Big[. $$

\Proof. First [EJ90, 6.6] implies that 
$l(\lambda^0)=\{ u_0, \ldots, u_m\} \dot\cup ]u_m,\infty[$
with $m = r - k-1,$ so that it remains to determine $u_0$. 

Let $u := u_0$. To compute this parameter, we recall from [EJ90] that 
$\le(\lambda) = 1$, and that the highest weight of 
$J(\lambda)$ is given by the unique dominant integral functional 
$\mu_0$ conjugate to $\lambda - \gamma_r$ 
under the Weyl group of $\hat\Delta_k$. 
Since in our context $\lambda$ and $\gamma_r$ can be identified 
with elements of $\b^*$ and the Weyl group contains in particular 
the permutations $\sigma_{ij}$ for two element $\gamma_i$ and $\gamma_j$, 
it follows that 
$$ \mu_0 
= \sigma_{r,k+1}.(\lambda - \gamma_r) 
= \lambda - \gamma_{k+1}. $$
That this functional is dominant is a consequence of the fact that 
$m_{k+1} - 2 \geq m_k$ 
which follows from (5.2). 

To calculate $u_0$, we use the fact that $L(\lambda^{u_0})$ and 
$L(\mu_0)$ correspond to the same eigenvalues of the Casimir operator 
of $\g$, which implies that 
$$ \| \lambda + \hat\rho\|^2 = \|\mu_0 + \hat\rho\|^2 = \|\lambda + \hat\rho - \gamma_{k+1}
\|^2$$
for $\hat\rho\:={1\over 2}\sum_{\hat\alpha\in\hat\Delta^+} \hat\alpha$. 
This means that 
$2 \la \lambda + \hat\rho, \gamma_{k+1} \ra=  \|\gamma_{k+1}\|^2$
which is the same as 
$$ (\lambda + \hat\rho)(\check \gamma_{k+1}) = 1. $$
Since $\lambda(\check \gamma_{k+1}) = \lambda(\check \gamma_r) = -u_0$, 
we obtain 
$u_0 = \hat\rho(\check \gamma_{k+1}) -1.$

To compute this number, we first observe that it equals one half of the 
trace of $\check\gamma_{k+1}$ on the positive root spaces in $\g_\C$. 
Since $\check\gamma_{k+1}$ corresponds to the element 
$2L(c_{k+1})$ under the Cayley transform, the contribution of the 
non-compact root spaces is given by 
$$ \tr_V L(c_{k+1}) = (r-1){d\over 2} + 1 $$
(this is a consequence of the Peirce decomposition of $V$; see [FaKo94, Th.\ IV.2.1]). 
The contribution of the compact roots equals 
$$ {1\over 2} \sum_{\alpha \in \Sigma^+} (\dim \h^\alpha) \alpha(2L(c_{k+1}))
= {d\over 2} \sum_{i < j} (\eps_j - \eps_i)(L(c_{k+1}))
= {d\over 2}(k - (r-k-1)) = (2k-r+1){d\over 2}.$$
Summing up, we obtain 
$$ u_0 = \hat\rho(\check \gamma_{k+1}) -1 
= (r-1 + 2k - r +1){d\over 2} = kd. 
\qeddis 

\subheadline{Support properties}

\Lemma V.6. Let $D \in {\cal S}'(V)$ be a tempered distribution whose support is contained 
in $\oline\Omega$. Let further $I\subeq \Pol(V)$ be an ideal such that 
$\la D, e^{-\tr} I \ra = \{0\}.$
Then 
$$ \supp(D) \subeq \{ x \in V \: (\forall f \in I) f(x) = 0\}. $$

\Proof. Since $D$ is supported by $\oline\Omega$, the 
product $D e^{-\tr}$ 
defines a tempered distribution on $V$ (cf.\ [HiNe99, Lemma B.1]).
Let $f \in I$. We claim that $f D= 0$. For all polynomials 
$p \in \Pol(V)$ we have 
$$ \la f D, e^{-\tr} p \ra = \la D, e^{-\tr} p f \ra = 0. $$
On the other hand, the Fourier transform ${\cal F}(fD) \in \Hol(T_\Omega)$ 
(cf.\ Section II) satisfies 
$$ \Big(p\big({\partial \over \partial z}\big) {\cal F}(fD)\Big) (ie) 
=\la fD, e^{-\tr} p(i\cdot) \ra = 0. $$
Since ${\cal F}(fD)$ is a holomorphic function on $T_\Omega$ and $p$ was arbitrary, 
we conclude that ${\cal F}(fD) = 0$, and hence that $fD = 0$. This implies that 
$\supp(D) \subeq f^{-1}(0)$ ([HiNe99, Lemma B.4]), hence yields the lemma. 
\qed

\Lemma V.7. Let $D \in {\cal S}'(V) \otimes F(\lambda)$ be a 
tempered distribution whose support is contained 
in $\oline\Omega$. If 
$$ \la D, e^{-\tr} I_k \otimes F(\lambda) \ra = \{0\}, $$
then $\supp(D) \subeq {\cal O}_k$. 

\Proof. Let $\alpha \in F(\lambda)^*$. Our assumption implies that 
$\la \alpha \circ D, e^{-\tr} I_k\ra = \{0\},$
so that Lemma V.6 implies that 
then $\supp(\alpha \circ D) \subeq {\cal O}_k$. Since $\alpha$ was arbitrary, the 
assertion follows. 
\qed

\Proposition V.8. If $L(\lambda)$ is unitarizable and spherical, then 
$\supp(R_\lambda^H) \subeq \supp(R_\lambda).$

\Proof. If $\supp(R_\lambda) = \oline\Omega$, there is nothing to show. 
So we may assume that $u = k{d\over 2}$ and 
$\supp(R_\lambda) = \oline{{\cal O}_k}$ (Theorem IV.4). Then 
$$e^{-\tr} I_k \otimes F(\lambda)\subeq J(\lambda), $$
so that the assumption that $L(\lambda)$ is spherical implies that 
$R_\lambda^H$ vanishes on this space. Hence the assertion follows from 
Lemma V.7. 
\qed

As the preceding proposition suggests, information on the support of 
$R_\lambda^H$ might provide information on whether $L(\lambda)$ is spherical or not. 

\Proposition V.9. About the support of the distributions $R_\lambda^H$, 
$\lambda = \lambda^s$, we have the following information: 
\item{(i)} If $s = jd$, $j \in \{ k, \ldots, r-1\}$, then
$R_\lambda^H$ is a  measure supported by $\oline{{\cal O}_k}$ for which 
$\oline{{\cal O}_{k-1}}$ is zero set. 
\item{(ii)} If $s > (r-1)d$, then $R_\lambda^H$ is a measure supported by 
$\oline\Omega$ for which $\partial \Omega$ is a zero set. 
\item{(iii)} If  $s\not\in\{0,\ldots,(r-1)d\} - 2\N_0$ 
and $s\leq(r-1)d$, 
then $R_\lambda^H$ is not a measure and its support is all of $\oline\Omega$.
\item{(iv)} If  $s = j d - 2 m$ with $m \in \N$ and 
$j \in \{k, \ldots, r -1 \}$, then $\supp(R_\lambda^H) \subeq \oline{{\cal O}_j}$. 

\Proof. (i) follows directly from Theorem V.1 and [HiNe99, Prop.\ III.6]. 

\nin(ii) This is a consequence of the fact that the density of $R_\lambda^H$ 
with respect to the $H$-semiinvariant measure $R_{s\over 2}$ on $\Omega$
does not vanish anywhere.

\nin(iii) In view of the relation 
$$ \Delta^l R_s = (s)_l R_{s + l}, \quad s \in \C, l \in \N $$
([HiNe99, Prop.\ III.1]), we obtain for $s > (r-1)d$ the relation 
$$ \eqalign{ \Delta^l R_{\lambda^s}^H 
&= \big({\textstyle s \over 2}\big)^{-1}_{-{\bf m}} \Delta^l \eta R_{s\over 2} 
= \big({\textstyle s \over 2}\big)^{-1}_{-{\bf m}} ({\textstyle s \over 2})_l \eta R_{s+2l\over 2} \cr
&= \big({\textstyle s \over 2}\big)^{-1}_{-{\bf m}} ({\textstyle s \over 2})_l 
\big({\textstyle s +2l\over 2}\big)_{-{\bf m}} R_{\lambda^{s + 2l}}^H 
= \big({\textstyle s \over 2} - {\bf m}\big)_l R_{\lambda^{s + 2l}}^H. \cr} $$

Since both sides are weakly holomorphic functions with values in 
${\cal S}'(V,F_0)$, we conclude that 
$$ \Delta^l R_{\lambda^s}^H = ({\textstyle s \over 2})_l R_{\lambda^{s + 2l}}^H \leqno(5.3) $$
holds for all $s \in \C$. 
Now the hypotheses shows that $({\textstyle s \over 2})_l\not=0$ for any $l\in \N_0$. 
Therefore (i) and (ii) together with (5.3) imply that 
$\oline\Omega \subeq \supp(R_{\lambda^s}^H)$, hence equality. 

\nin (iv) From the proof of [HiNe99, Prop.\ V.2] we get 
$$ \Delta({\partial \over \partial x}) R_{\lambda^s}^H = R_{\lambda^{s-2}}^H, $$
so that (iv) follows from (i). 
\qed

The situation is particularly simple for $\lambda^0 = 0$, where $k = 0$, so that 
Proposition V.9 gives complete information on the support of $R_{\lambda^s}^H = R_{s\over 2}$.
Fortunately the classification in Theorem V.5 shows that the information from 
Proposition V.9 suffices to deal with all those 
distributions $R_\lambda^H$ for which $F(\lambda)$ is spherical and $L(\lambda)$ 
is unitary. This is made precise by the following result: 

\Proposition V.10. Suppose that $F(\lambda^0)$ is $\h \cap \k$-spherical and 
that $L(\lambda)$ is singular and unitary. Then the following assertions 
hold: 
\item{(i)} $R_\lambda^H$ is a measure if and only if $\le(\lambda)$ is odd. 
\item{(ii)} If $u = j{d\over 2}$ with $j \leq r-1$ odd, then 
$R_\lambda^H$ does not vanish on $e^{-\tr} I_j \otimes F(\lambda)$. 

\Proof. (i) Let $\lambda = \lambda^{u_j}$. Then $\le(\lambda) = j+1$ and 
${u_j\over 2} = \big(k + {j \over 2}\big){d\over 2}$ (Theorem V.5). 
Since 
$$ {u_j\over 2} \leq \big(k + {r-k-1 \over 2}\big){d\over 2}
= {r+k-1 \over 2}{d\over 2}
\leq (r-1){d\over 2},  $$
the corollary follows from Proposition V.9(i). 

\par\nin (ii) The ideal $I_j$ contains in particular the polynomial 
$\Delta_{{\bf j+1}}$ for 
${\bf j+1}=(\underbrace{1,\ldots,1}_{(j+1)-times}, 0,\ldots,0)$, so that 
it suffices to show that 
$$ \la R_\lambda^H, e^{-\tr} \Delta_{{\bf j+1}} v_\lambda \ra \not =0, $$
where $v_\lambda \in F(\lambda)$ is a real highest weight vector for 
$\h$ with highest weight 
$\lambda_\b = \sum_{i=1}^r {m_i - u\over 2} \eps_i.$

Let ${\bf s} := (s_1, \ldots, s_r) \in \C^r$ and $R_{\bf s} \in 
{\cal S}_{\oline\Omega}(\C)$ the multiparameter Riesz distribution which is 
uniquely determined by 
$$ {\cal F}(R_{\bf s})(ix) = \Delta_{\bf s}(x^{-1}) \quad \hbox{ for } \quad 
x \in \Omega $$
(cf.\ [FaKo94, Ch.\ VII]). Now [HiNe99, Lemma V.4] yields 
$$ \la R_\lambda^H, v_\lambda \ra = \la v_0, v_\lambda \ra 
R_{{u - {\bf m}\over 2}}. $$
On the other hand we have 
$$ \Delta_{{\bf j+1}}.R_{{u - {\bf m}\over 2}}
= ({u - {\bf m}\over 2})_{{\bf j+1}} R_{{u - {\bf m}\over 2}+ {\bf j+1}} $$
with 
$$ ({u - {\bf m}\over 2})_{{\bf j+1}} 
=  \prod_{l = 1}^{j+1} \Big({u-m_l\over 2} - (l-1){d\over 2}\Big)$$
([HiNe99, Prop.\ III.1(iii)]). This eventually leads to 
$$ \eqalign{ \la R_\lambda^H, e^{-\tr} \Delta_{{\bf j+1}} v_\lambda \ra 
&= {\cal F}(\Delta_{{\bf j+1}} \la R_\lambda^H, v_\lambda \ra)(-ie) 
=  ({u - {\bf m}\over 2})_{{\bf j+1}} 
\Delta_{{u - {\bf m}\over 2}+ {\bf j+1}}(e) \la v_0, v_\lambda \ra \cr
&=  ({u - {\bf m}\over 2})_{{\bf j+1}} \la v_0, v_\lambda \ra. \cr} $$
Since $\la v_0, v_\lambda \ra \not= 0$ (cf.\ [HiNe99, Lemma II.1]), 
it remains to show that 
$({u - {\bf m}\over 2})_{{\bf j+1}}$
does not vanish if $j$ is odd. In view of $m_l = 0$ for $l > k$, we have 
$$ \prod_{l = k+1}^{j+1} \Big({u-m_l \over 2} - (l-1){d\over 2}\Big)
= \prod_{l = k+1}^{j+1} {d\over 2} \big({j \over 2} - (l-1)\big)\not= 0. $$
For $l \leq k$ we use $j \geq 2k$ (Theorem V.5) to obtain 
$$ -m_l + {u \over 2} -(l-1){d\over 2} 
\geq {u \over 2} -(l-1){d\over 2} 
\geq {d\over 2}(k -(l-1)) > 0. $$
Therefore all factors in 
$({u - {\bf m}\over 2})_{{\bf j+1}}$ are non-zero. 
\qed

\Theorem V.11. {\rm(Classification for singular supports)} Let $G$ be a 
simply connected hermitian Lie group of tube type and $G/H$ an associated 
symmetric space of Cayley type. 
Suppose that $F(\lambda)$ is $\h \cap \k$-spherical and 
that $L(\lambda)$ is unitary with $\supp(R_\lambda) \not=\oline\Omega$. 
Then $L(\lambda)$ is spherical if and only if 
$\le(\lambda)$ is odd. 

\Proof. Since $\supp(R_\lambda)$ is smaller than $\oline\Omega$, we have 
$u = j{d\over 2}$ for some 
$j \leq r-1$ (Theorem IV.4). Moreover, in view of 
Remark II.3, Theorem V.5 shows that $\le(\lambda)$ is odd if and only if 
$j$ is even. If this is the case, then $L(\lambda)$ is spherical by 
Theorem II.8. If $j$ is odd, then Proposition V.10 shows that 
$R_\lambda^H$ does not vanish on $e^{-\tr} I_j \otimes F(\lambda)$ 
which is contained in $J(\lambda)$ (Proposition IV.6). Therefore Corollary V.4 implies that 
$L(\lambda)$ is not spherical. 
\qed

\Corollary V.12. {\rm(Classification for the scalar case)} 
If $L(\lambda)$ is singular, unitary and of scalar type, then 
$L(\lambda)$ is spherical if and only if $\le(\lambda)$ is odd. 

\Proof. In the scalar case $F(\lambda)$ is always $H\cap K$-spherical and 
$\lambda$ is singular if and only $\supp(R_\lambda) \not=\oline\Omega$. 
Therefore the assertion follows from Theorem V.11. 
\qed

\Lemma V.13. If $L(\lambda)$ and $L(\mu)$ are unitary and spherical, then the same holds 
for $L(\lambda + \mu)$. 

\Proof. Since $L(\lambda)$ and $L(\mu)$ are unitary, the highest weight module 
$L(\lambda + \mu)$ can be embedded into the unitary module 
$L(\lambda) \otimes L(\mu)$ 
as the submodule generated by 
$v_\lambda \otimes v_\mu$, where $v_\lambda$ and $v_\mu$ are highest weight 
vectors ([Ne99, Prop.\ IX.1.15]). It follows in particular that 
$L(\lambda + \mu)$ is unitary. Moreover, if 
$\nu_\lambda \in \big(L(\lambda)^\sharp\big)^\h$ and 
$\nu_\mu \in \big(L(\mu)^\sharp\big)^\h$ are non-trivial 
$\h$-invariant antilinear functionals, then 
$\nu_\lambda \otimes \nu_\mu$ is an $\h$-invariant antilinear functional 
on $L(\lambda) \otimes L(\mu)$ with 
$$ (\nu_\lambda \otimes \nu_\mu)(v_\lambda \otimes v_\mu) 
= \nu_\lambda(v_\lambda) \nu_\mu(v_\mu) \not= 0 $$
(Proposition I.5). This shows that the restriction of 
$\nu_\lambda \otimes \nu_\mu$ to $L(\lambda + \mu)$ is non-zero and hence that  
this module is spherical. 
\qed

\Remark V.14. Let $u_0 = kd$ as in Theorem V.5. Then Theorem V.11 shows that 
$L(\lambda)$ is spherical for $u = u_j = kd + j{d\over 2}$ if $j$ is even 
(which corresponds to odd reduction level) and not spherical if 
$j$ is odd, provided $u_j \leq (r-1){d\over 2}$. This means that 
$$  j \leq r-1 - 2k. $$
For an odd $j \in \{r-2k, \ldots, r - k - 1\}$ Theorem V.11 gives no information.

Suppose that for some $j$ the module $L(\lambda^{u_j})$ is spherical. 
We know from Corollary V.12 that $L(d\zeta)$ is spherical, so that 
Lemma V.13 and 
$$ \lambda^{u_{j+2}} = \lambda^{u_j} + d\zeta $$
imply that $L(\lambda^{u_{j+2}})$ is spherical. 
This argument can be iterated, and it shows that 
if for one odd value of $j$ the module $L(\lambda^{u_j})$ is spherical, then 
it is true for all larger (odd) values of $j$. 
\qed

We have already seen in Corollar V.12 that for singular support of $R_\lambda$ 
and even reduction level the scalar type highest weight modules are not spherical. 
It is instructive to observe that this in general not ruled out by the support 
condition in Proposition V.8. 
The following proposition shows that the support behavior is actually quite 
complicated. 

\Proposition V.15. For $s > (r-1){d\over 2}$ we always have 
$\supp(R_{s\over 2}) \subeq \supp(R_s).$ 
For $s = {kd\over 2}$, $k \in \{0, 1,\ldots, r-1\}$ we have: 
\item{(i)} If $k$ is even, then 
$\supp(R_{s\over 2}) \subeq \supp(R_s).$
\item{(ii)} If $k$ is odd, then we have: 
\itemitem{(a)} If $d = 1,2$, then 
$\supp(R_{s\over 2})\not\subeq \supp(R_s)$.
\itemitem{(b)} If $d = 4$, then 
$\supp(R_{{s\over 2}})\subeq \supp(R_{s})$.
\itemitem{(c)} If $d =8$ and $r = 3$, then 
$\supp(R_{{s\over 2}}) \subeq \supp(R_{s}).$
\itemitem{(d)} If $d =n-2$ and $r = 2$, then 
$\supp(R_{{s\over 2}}) \subeq \supp(R_{s})$ 
if and only if $n \equiv 2 \mod 4$. 

\Proof. (i) This follows from 
$\supp(R_{{s\over 2}}) = \oline{{\cal O}_{k\over 2}} 
\subeq \oline{{\cal O}_{k}} = \supp(R_s).$

\par\nin (ii) (a) $d = 1$: Then 
$$ {kd \over 4} = {k \over 4}
\not\in \{j {\textstyle{d\over 2}} \: j =0,\ldots, r-1\} - \N_0 $$
because $k \not\equiv 2j \mod 4$ for each $j$. Hence 
Proposition V.8 implies that 
$$ \supp(R_{{s\over 2}}) = \oline{\Omega} \not\subeq \supp(R_{s}) $$
in this case. 

\par\nin $d =2$: If $k$ is odd, then 
$$ {kd \over 4} = {k \over 2}
\not\in \{j {\textstyle{d\over 2}} \: j =0,\ldots, r-1\} - \N_0 
= \{j \: j =0,\ldots, r-1\} - \N_0 $$
because ${k \over 2}$ is not an integer. 

\par\nin (iii) $d =4$: If $k$ is odd, then 
$$ {kd \over 4} = k 
\in \{j {\textstyle{d\over 2}} \: j =0,\ldots, r-1\} - \N_0 
= \{2 j \: j =0,\ldots, r-1\} - \N_0 $$
and the minimal $j$ with 
$k \in 2 j - \N_0$ is given by 
$j = {k + 1\over 2} \leq k$. We conclude that in this case we always have 
$\supp(R_{{s\over 2}}) \subeq \oline{{\cal O}_j} 
\subeq \oline{{\cal O}_k} = \supp(R_{s}).$

\par\nin (iv) $d =8$ and $r = 3$: Then $k$ odd implies that $k =
1$. Now 
$$ {kd \over 4} = 2 = 4 \cdot 1 - 2 
\in \{j {\textstyle{d\over 2}} \: j =0,\ldots, r-1\} - \N_0 
= \{4 j \: j =0,1,2\} - \N_0. $$
Hence 
$ \supp(R_{{s\over 2}}) \subeq \oline{{\cal O}_1} 
=  \supp(R_{s}).$

\par\nin (v) $d =n-2$ and $r = 2$: Then $k = 1$ and therefore 
$$ {kd \over 4} = {n - 2 \over 4} 
\in \{j {\textstyle{d\over 2}} \: j =0,\ldots, r-1\} - \N_0
= \{j {\textstyle{d\over 2}} \: j =0,1\} - \N_0 
= \{0, {n-2\over 2}\} - \N_0 $$
is equivalent to 
$$ {n - 2 \over 2} - {n - 2 \over 4} 
= {n - 2 \over 4} \in \N. $$
Thus 
$\supp(R_{{s\over 2}}) \not\subeq \supp(R_{s})$
for $n \not\equiv 2 \mod 4$, and for 
$n \equiv 2 \mod 4$ we have 
$\supp(R_{{s\over 2}}) \subeq \oline{{\cal O}_1} 
=  \supp(R_{s}).$ 
\qed

\Example V.16. We consider the Lie algebra $\g= \su(r,r)$ corresponding to the Jordan algebra 
$V = \Herm(r,\C)$. Here $\h \cong \R \oplus \sL(r,\C)$, which we consider as a quotient 
of $\gl(r,\C)$ modulo $i\R$. This is most natural because $\gl(r,\C)$ acts naturally on 
$V$ by $X.A = XA + AX^*$, and the kernel of this representation is $i\R$. 
The complexification of $\h$ is 
$$ \h_\C \cong \C \oplus \sL(r,\C) \oplus \sL(r,\C) $$
which we consider as a quotient of $\gl(r,\C) \oplus \gl(r,\C)$ acting on 
$V_\C \cong \End(\C^r)$ by $(X_1, X_2).A = X_1 A + A X_2^*$. 
Accordingly we write the highest weights of simple $\h_\C$ modules as 
$$ \lambda = (\lambda_1, \ldots, \lambda_{2r}) = (\lambda^+, \lambda^-) $$
with respect to the positive system 
$$ \hat\Delta_k^+ = \{ \eps_i - \eps_j \: 1 \leq i < j \leq r; r+1 \leq i < j \leq 2r \} $$
of roots of $\sL(r,\C) \oplus \sL(r,\C)$. 

The corresponding representation space $F(\lambda)$ has the structure 
$$ F(\lambda) = F(\lambda^+) \otimes F(\lambda^-) 
= F(\lambda^+) \otimes F(\tilde\lambda^-)^*
\cong \Hom(F(\tilde \lambda^-), F(\lambda^+)), $$
where $\tilde\lambda^- = (-\lambda_{2r}, \ldots, -\lambda_{r+1})$. 
The action of $\h$ on this space is given by 
$$ X.A = \pi_{\lambda^+}(X) A + A \pi_{\tilde\lambda^-}(X)^*. $$
Therefore a fixed vector for $\h \cap \k \cong \su(r,\C)$ corresponds to an intertwining 
operator for $\su(r,\C)$, and hence for $\sL(r,\C)$, in 
$\Hom(F(\tilde \lambda^-), F(\lambda^+))$. This means that 
$F(\lambda)$ is spherical if and only if 
$F(\lambda^+) \cong F(\tilde \lambda^-)$ as $\sL(r,\C)$-modules. 
With $\mu := \lambda^+$ we therefore get 
$F(\lambda) \cong \End(F(\mu))$ with $\sL(r,\C)$ acting by 
$X.A = \pi_{\mu}(X) A + A \pi_{\mu}(X)^*.$

The restricted highest weight 
of this representation with respect to the subspace 
$\b$ of real diagonal elements in $\h$ is given by 
$$ \lambda_\b = \sum_{j=1}^r {m_j \over 2} \eps_j = 
\sum_{j=1}^r {\lambda_{r+1-j} - \lambda_{r+j} \over 2}\eps_j. $$
The normalization $m_r = 0$ 
leads to $\lambda_1 = \lambda_{2r}$, and for the spherical representations further to 
$$ \lambda_\b = \sum_{j=1}^r \lambda_{r+1-j} \eps_j $$
because $\lambda^+ = - \tilde\lambda^-$ in this case. 
We determine $k \in \{0,\ldots, r-1\}$ by 
$ \lambda_1 = \lambda_2 = \ldots = \lambda_{r-k} > \lambda_{r-k+1}. $
Then Theorem V.5 implies that 
$m = r-1-k$ and $u_0 = 2k$ because $d = 2$ holds for $\g = \su(r,r)$. 

It is interesting to observe that the relation 
$F(\lambda) \cong \End(F(\mu))$ as representations of $H$ implies that 
$\tau_\mu \circ \tilde P \: T_\Omega \to \End(F(\mu))$ 
is $H$-equivariant, so that the injectivity of the Fourier transform implies 
for $v_0 = \1 \in \End(F(\mu))$ the relation 
$$R_\mu = R_\lambda^H, $$
where $F(\mu)$ is interpreted as a module of $\h_\C$ which is trivial on the second 
$\sL(r,\C)$-factor. It has been shown by J.\ L.\ Clerc in [Cl95] that 
for these representations $L(\mu)$ is unitary if and only if 
${u\over 2} \in \{k,k+1, \ldots, r-1\}$ or $u > (r-1)$. According to Proposition V.8, these 
are the cases where $R_\mu = R_\lambda^H$ is a measure, and they correspond to the 
cases with odd reduction level. Then the unitary $\g$-module  
$L(\mu) \otimes L(\mu)^*\subeq \End(L(\mu))$ contains a $\g$-submodule isomorphic 
to $L(\lambda)$, and the trace is an $\h$-invariant linear functional on $L(\lambda)$. 
This construction can also be used to prove that the representations $L(\lambda)$ 
are spherical for odd reduction levels. 
\qed

Our results in this section suggest that for even reduction level a 
singular highest weight representation is never spherical. At least we do not know 
of any counterexample. 

Similarly, in all cases of singular representation where we could decide whether $L(\lambda)$ 
is spherical or not, this happens if and only if the distribution 
$R_\lambda^H$ is a measure. In the scalar case this follows from Proposition V.9 and Theorem V.5. 

\subheadline {Appendix: Regularity of the holomorphic discrete series}

In this appendix we show that the highest weight modules $L(\lambda)$ corresponding to 
the discrete series of a compactly causal symmetric space $G/H$ are regular in the 
sense that $N(\lambda) \cong L(\lambda)$. This shows in particular that these 
representations do not provide any further information on sphericalness of singular 
highest weight representations. 

Let $\g$ be a hermitian Lie algebra, $\t\subeq \g$ a compactly embedded Cartan subalgebra 
and $\hat\Delta=\hat\Delta(\g_\C,\t_\C)$ the corresponding root system. 
We choose a positive system $\hat\Delta^+$ of $\hat\Delta$ as in Section I and set 
$\hat\rho\:={1\over 2}\sum_{\hat\alpha\in\hat\Delta^+} \hat\alpha$, 
$\rho_k\:={1\over 2}\sum_{\hat\alpha\in\hat\Delta_k^+} \hat\alpha$ and 
$\hat\rho_n\:={1\over 2}\sum_{\hat\alpha\in\Delta_n^+} \hat\alpha$.

If $\lambda\in i\t^*$ is dominant integral w.r.t. $\hat\Delta_k^+$, then the condition 
for $L(\lambda)$ to correspond to the relative holomorphic discrete series of $G$ is 
given by Harish Chandra's condition (cf.\ [HC56])
$$(\forall\hat\alpha\in \hat\Delta_n^+) \qquad \la \lambda+\hat\rho,\hat\alpha\ra<0.\leqno DS(G)$$
Elementary algebraic considerations involving the Parthasarathy inequality (cf.\ [EHW83, Prop.\ 3.9])
or [Jan79,p.\ 35] imply that DS(G) implies that $\lambda$ is regular, i.e, $N(\lambda)
=L(\lambda)$.  

Equip now  $\g$ with an involution $\tau$ such that $(\g,\tau)$ becomes compactly causal. 
Recall the root system $\Delta=\Delta(\g^c,\a)$. The fact that $\z(\k)\subeq\q$ implies that 
we can choose $\hat\Delta^+$ such that $-\tau\hat\Delta^+\subeq \hat\Delta^+\cup\a^\bot$. 
Then $\Delta^+\:=\hat\Delta^+\res_\a\bs\{0\}$ 
is a positive system of $\Delta$. Recall the definition of 
$\Delta_n^+$ and $\Delta_k^+$ from Section I. 
Set $\rho\:={1\over 2}\sum_{\alpha\in\Delta^+} m_\alpha \alpha$, 
$\rho_k\:={1\over 2}\sum_{\alpha\in\Delta_k^+} m_\alpha \alpha$ and 
$\rho_n\:={1\over 2}\sum_{\alpha\in\Delta_n^+} m_\alpha \alpha$ with $m_\alpha=\dim(\g^c)^\alpha$. 

\par Let $p\: \t_\C\to \a_\C$ denote the orthogonal projection with respect to the 
Cartan--Killing form. Then the adjoint map $p^*\: \a_\C^*\to \t_\C^*$ is injective 
and via this inclusion mapping we identify in the sequel $\a_\C^*$ with a subspace 
of $\t_\C^*$.

\par For a $H\cap K$-spherical highest weight $\lambda\in \a^*\subeq i\t^*$ the condition for $L(\lambda)$
to belong to the relative holomorphic discrete series of $G/H$ is given by 

$$(\forall\alpha\in \Delta_n^+) \qquad \la \lambda+\rho,\alpha\ra<0\leqno DS(G/H)$$
(cf.\ [\'O\O91], [H\'O\O91]; see also [Kr99b]). 
It follows easily from the fact that $\hat\Delta^+$ and $\Delta^+$ are compatible that 
$DS(G/H)$ implies $DS(G)$. Conversely it was observed in [\'O\O{}88, Lemma 7.4] that 
$DS(G)$ is a weaker condition than $DS(G/H)$ for $H\cap K$-spherical highest weights $\lambda$.
Hence the following result is of interest.

\Theorem A.1. Let $\lambda\in \a^*$ be dominant integral for 
$\hat\Delta_k^+$ and assume that $DS(G/H)$ holds. Then $\lambda$ is regular, i.e., 
$L(\lambda)=N(\lambda)$.\qed 

We are going to give a simple algebraic proof of Theorem A.1 using the Parthasarathy inequality.

\Lemma A.2. We have $\hat\rho_n=\rho_n$. 

\Proof. This follows from $\hat\rho_n\in i\z(\k)^*$, $\z(\k)\subeq \a$ and 
$\hat\Delta_n^+\res_\a \subeq \Delta^+$. \qed

\Lemma A.3. If $\lambda\in \a^*$ and $\mu\in i\t^*$, then 
$$\|\lambda+\hat\rho\|^2-\|\mu+\hat\rho\|^2
= \|\lambda+\rho\|^2-\|\mu+\rho\|^2
- 2\la \hat\rho_k-\rho_k,\mu\ra.$$

\Proof. By Lemma A.2 we have 
$$\eqalign{
&\ \ \ \  \|\lambda+\hat\rho\|^2-\|\mu+\hat\rho\|^2
=\|\lambda+\rho+(\hat\rho-\rho)\|^2-\|\mu+\rho+(\hat\rho-\rho)\|^2\cr
&=\|\lambda+\rho+(\hat\rho_k-\rho_k)\|^2-\|\mu+\rho+(\hat\rho_k-\rho_k)\|^2\cr
&=\|\lambda+\rho\|^2 +2\la \lambda+\rho, \hat\rho_k-\rho_k\ra 
-\|\mu+\rho\|^2- 2\la \mu+\rho, \hat\rho_k-\rho_k\ra\cr
&=\|\lambda+\rho\|^2 -\|\mu+\rho\|^2- 2\la \mu, \hat\rho_k-\rho_k\ra,\cr}$$
where the last equality comes from $\hat\rho_k-\rho_k\in \a^\bot$. \qed

\Lemma A.4.   We have 
$\hat\rho_k=\rho_k+\rho_k^0$
with $\rho_k^0\:={1\over 2}\sum_{\alpha\in \hat\Delta_k^+\cap\a^\bot}\alpha\in \a^\bot$. 

\Proof. Let $\alpha\in \hat\Delta_k^+$ and rewrite 
$\alpha$ as $\alpha=\alpha_\q+\alpha_\h$ with $\alpha_\q=\alpha\res_{\t\cap\q}$ and 
$\alpha_\h=\alpha\res_{\t\cap\h}$. Assume that 
$\alpha_\q\neq 0$. Then 
$(-\tau).\hat\Delta_k^+ \subeq \hat\Delta_k^+ \cup \a^\bot$
and the fact that 
$\hat\Delta_k^+$ and $\Delta_k^+$ are compatible yield 
$-\tau.\alpha=\alpha_\q -\alpha_\h\in \hat \Delta_k^+.$ 
Hence we get 
$${1\over 2}\sum_{\alpha\in \hat\Delta_k^+\atop \alpha_\q\neq 0}\alpha=\rho_k$$ 
Now the assertion of the lemma follows from 
$\hat\rho_k={1\over 2}\sum_{\alpha\in \hat\Delta_k^+\atop \alpha_\q\neq 0}\alpha+
\rho_k^0.$ \qed

\nin {\bf Proof of Theorem A.1.} In view of Parthasarathy's condition (cf.\ [EHW83, Prop.\ 3.9]), 
we only have to check that 
$$\|\lambda+\hat\rho\|^2-\|\mu+\hat\rho\|^2<0$$
holds for all $\k_\C$-highest weights $\mu\neq \lambda$ of $N(\lambda)$. 
Note that $\mu$ can be written as 
$\mu=\lambda-\sum_{\alpha\in \hat\Delta_n^+} k_\alpha \alpha$ with 
$k_\alpha\in \N_0$. By Lemma A.2 we then have have 
$$\eqalign{\|\lambda+\hat\rho\|^2-\|\mu+\hat\rho\|^2&= 
\|\lambda+\rho\|^2-\|\mu+\rho\|^2
- 2\la \hat\rho_k-\rho_k,\mu\ra\cr 
&=-\|\sum_{\alpha\in \hat \Delta_n^+} k_\alpha \alpha\|^2 +
2\Big(\sum_{\alpha\in \hat\Delta_n^+} k_\alpha \la \lambda+\rho, \alpha\ra\Big) -
2\la \hat\rho_k-\rho_k,\mu\ra.\cr}$$
In the bottom line the first summand is clearly negative, 
the second one by assumption, and the non-positivity of the third one 
follows from Lemma A.4 since $\mu$ is dominant with respect to 
$\hat\Delta_k^+$. 
\qed

\def\entries{

\[Ba88 van den Ban, E., {\it The principal series for a reductive symmetric space. I. $H$-fixed distribution vectors}, Ann. Sci. \'Ecole Norm. Sup. (4) {\bf 21:3} (1988), 359--412

\[BaDe88 van den Ban, E., and P.\ Delorme, {\it Quelques propri\'et\'es des repr\'esentations
sph\'eriques pour les espaces sym\'etriques r\'eductifs}, J. Funct. Anal. {\bf 80} (1988), 
284--307

\[Bi90 Bien, F., ``${\cal  D}$-modules and spherical representations,''
Mathematical Notes {\bf 39}, Princeton University Press, Princeton, NJ, 1990

\[BrDe92 Brylinski, J.L., and P. Delorme, {\it Vecteurs distributions H-Invariants pur 
les s\'eries principales g\'en\'eralis\'ees d'espaces  sym\'etriques reductifs et 
prolongement meromorphe d'int\'egrales d'Eisenstein}, Invent. math. {\bf 109} (1992), 619--664

\[Cl95 Clerc, J.--L., {\it Laplace transform and unitary highest weight modules}, 
J.\ Lie Theory {\bf 5} (1995), 225--240

\[DES91 M.\ G.\ Davidson, T.\ J.\ Enright, and R.\ J.\ Stanke, ``Differential 
operators and highest weight representations," Memoirs of the 
Amer.\ Math.\ Soc., Providence, 1991 

\[DG93 Ding, H., and K.\ Gross, {\it Operator-valued Bessel functions on Jordan 
algebras}, J. reine ang. Math. {\bf 435} (1993), 157--196

\[EJ90 Enright, T. J., and  A. Joseph, {\it An intrinsic analysis of unitarizable highest 
weight modules}, Math. Ann. {\bf 288} (1990), 571--594

\[EHW83 Enright, T. J., R. Howe, and N. Wallach, {\it A classification of 
unitary highest weight modules} in ``Representation Theory of Reductive 
Groups,'' Park City, UT, 1982, pp. 97--149; Progr. Math. {\bf 40} (1983), 
97--143 

\[FaKo90 Faraut, J., and A.\ Koranyi, {\it Function Spaces 
and Reproducing Kernels on Bounded Symmetric Domains}, 
J.\ Funct.\ Anal. {\bf 88:1} (1990), 64--89

\[FaKo94 ---, ``Analysis on Symmetric Cones," 
Oxford Mathematical Monographs, Oxford University Press, 1994

\[FJ80 Flensted-Jensen, M., {\it Discrete series for semisimple
symmetric spaces}, Ann. of Math. {\bf 111:2} (1980), 253--311

\[HC56 Harish-Chandra, {\it Representations of semi-simple Lie groups,  VI}, 
Amer. J. Math. {\bf 78} (1956), 564--628

\[Hel78 Helgason, S., ``Differential Geometry, Lie Groups, and Symmetric 
Spa\-ces,'' Acad. Press, London, 1978

\[Hel84 ---, ``Groups and Geometric Analysis,'' Acad. Press, London,
1984

\[HiKr98 Hilgert, J., and B. Kr\"otz, {\it The Plancherel
Theorem for invariant Hilbert spaces}, Math.\ Z., to appear 

\[HiNe99 Hilgert, J., and K.--H. Neeb, {\it Vector valued Riesz distributions on Euclidean 
Jordan algebras}, J. Geom. Analysis, to appear 

\[Hi\'Ol96 Hilgert, J.\ and 
G.\ \'Olafsson, ``Causal Symmetric Spaces, Geometry and
Harmonic Analysis,'' Acad. Press, 1996 

\[H\'O\O{}91 Hilgert, J., G.\ Olafsson, and B. \O rsted, {\it Hardy Spaces on 
Affine Symmetric Spaces}, J. reine angew. Math. {\bf 415}(1991), 189--218

\[Jak83 Jakobsen, H.\ P., {\it Hermitian symmetric spaces and their unitary 
highest weight  modules}, J.\ Funct.\ Anal. {\bf 52} (1983), 385--412 

\[Jan79 Jantzen, J. C., ``Moduln mit einem h\"ochsten Gewicht," 
Lecture Notes in Math.\ {\bf 750}, Springer, 1979

\[KoWo65 Kor\'anyi, A., and J. A. Wolf, 
{\it Realization of hermitean symmetric 
spaces as generalized half planes}, Ann. of. Math. {\bf 81} (1965), 265--288

\[Kr99a Kr\"otz, B., {\it Norm estimates for unitary highest weight
modules}, Ann.\ Inst. Fourier {\bf 49:4)} (1999), 1242--1264

\[Kr99b ---, {\it Formal dimension for semisimple symmetric 
spaces}, Comp.\ Math., to appear 

\[KN\'O97 Kr\"otz, B., K. - H. Neeb, and G. \'Olafsson, {\it Spherical 
representations and mixed symmetric spaces}, Represent. Theory {\bf
1} (1997), 424--461

\[Lo69 Loos, O., ``Symmetric Spaces I : General Theory,'' 
Benjamin, New York, Amsterdam, 1969

\[MaOs84 Matsuki, T., and T. Oshima, {\it A description of discrete
series for semisimple symmetric spaces}, Adv. Stud. Pure Math. {\bf
4} (1984), 229--390 

\[Ne99  Neeb, K.--H., ``Holomorphy and Convexity in Lie Theory,''
Expositions in Mathematics {\bf 28}, de Gruyter, Berlin, 1999 

\[\'Ol87 \'Olafsson, G., {\it Fourier and Poisson transformation associated to a
semsisimple symmetric space}, Invent. math. {\bf 90} (1987), 605--629

\[\'O\O{}88  \'Olafsson, G., and B. \O{}rsted, 
{\it The holomorphic discrete series for affine symmetric spaces. I},  
J. Funct. Anal. {\bf 81:1} (1988), 126--159 

\[\'O\O91 ---, {\it The holomorphic discrete series of  
affine symmetric spaces and representations with reproducing kernels},
Trans. Amer. Math. Soc. {\bf 326} (1991), 385-405

\[VR76 Vergne, M., and H. Rossi, {\it Analytic continuation of the holomorphic 
discrete series of a semisimple Lie group}, Acta Math. {\bf 136} (1976), 1--59

\[Wal79 Wallach, N., {\it The analytic continuation of the discrete series I, 
II}, Transactions of the Amer. Math. Soc. {\bf 251} (1979), 1--17, 19--37

}

{\sectionheadline{\bf References}
\frenchspacing
\entries\par}

\dlastpage
\vfill\eject

\bye